\documentclass[reqno,12pt]{amsart}
\usepackage{amsthm}
\usepackage{amsmath, mathrsfs}
\usepackage{color}
\usepackage{amssymb}
\usepackage{hyperref}
\usepackage{enumerate}
\usepackage{latexsym,array}
\usepackage{amsfonts}
\usepackage{shadow}

\newcommand{\V}{\mathcal{V}}
\newcommand{\M}{\mathcal{M}}
\newcommand{\NN}{\mathcal{N}}
\newcommand{\LL}{\mathcal{L}}

\newtheorem{Pa}{Paper}[section]
\newtheorem{Tm}[Pa]{{\bf Theorem}}
\newtheorem{La}[Pa]{{\bf Lemma}}
\newtheorem{Dn}[Pa]{{\bf Definition}}
\newtheorem{Cy}[Pa]{{\bf Corollary}}

\newtheorem{Rk}[Pa]{{\bf Remark}}
\newtheorem{Pn}[Pa]{{\bf Proposition}}

\author[D. Alpay]{Daniel Alpay}
\address{(DA) Department of Mathematics\\
Ben–Gurion University of the Negev\\
Beer-Sheva 84105 Israel} \email{dany@math.bgu.ac.il}
\author[F. Colombo]{Fabrizio Colombo}
\address{(FC) Politecnico di
Milano\\Dipartimento di Matematica\\ Via E. Bonardi, 9\\20133
Milano, Italy}
\email{fabrizio.colombo@polimi.it}
\author[I. Sabadini]{Irene Sabadini}
\address{(IS) Politecnico di
Milano\\Dipartimento di Matematica\\ Via E. Bonardi, 9\\20133
Milano, Italy}
\email{irene.sabadini@polimi.it}
\title[Quaternionic Krein spaces]
{Inner product spaces and Krein spaces
in the quaternionic setting} \oddsidemargin 0.2in \evensidemargin
0.2in \topmargin -0.5in \textwidth 15.5truecm \textheight 23truecm
\def\H{\mathbb H}
\def\R{\mathbb R}

\def\(s){\mathscr S(\R\times\R)}

\def\W{\mathcal W}

 \keywords{Inner product spaces,
quaternionic topological vector spaces, Krein spaces, quaternionic
functional analysis.}


\thanks{{\sl Acknowlegments:} D. Alpay thanks the Earl Katz family for
endowing the chair
which supported his research, and the Binational Science
Foundation Grant number 2010117. F. Colombo and
I. Sabadini
acknowledge the Center for Advanced Studies of
the Mathematical Department of the Ben-Gurion University of the Negev for
the support and the kind hospitality during
the period in which part of this paper has been written. The authors thank V. Bolotnikov for useful discussions.}
\begin{document}
\maketitle \tableofcontents
\parindent 0cm
\begin{abstract}
In this paper we provide a study of quaternionic inner product
spaces. This includes ortho-complemented subspaces, fundamental decompositions
as well as a number of results of topological nature.
Our main purpose is to show that a closed uniformly  positive subspace in a quaternionic Krein space is
ortho-complemented, and this leads to our choice of the results presented in the paper.
\end{abstract}

\parindent 0cm

\section{Introduction}
\setcounter{equation}{0}
The purpose of this paper is to study quaternionic inner product spaces
and, in particular, Krein spaces. Quaternionic Hilbert spaces are known for a long time, see for
instance \cite{MR83k:46003}, and \cite{MR1333599,MR768240,MR0137500} for
various applications to quantum mechanics.

Some aspects of the theory of
quaternionic Pontryagin spaces have been studied in \cite{as3}.
The finite dimensional case is also of particular interest; see e.g.,
\cite{MR2001c:15021,MR97h:15020,MR2242728,MR2433159,MR2988215,MR2372587}.
Krein spaces are, roughly speaking, the
direct sum of two in general {\sl infinite dimensional} Hilbert spaces,
and in particular the previous references do not treat this case.
While preparing the work \cite{abcs2-milan} on interpolation of
Schur multipliers in the case of vector-valued
slice-hyperholomorphic functions, we realized that no reference
seemed to be available for a number of important results on
quaternionic Krein spaces. The motivation of the present paper
was to fill part of this gap. In the process, we found that we
needed to prove several
functional analysis results in the quaternionic setting.\\

In the complex case, the starting
point is a complex vector space $\V$, endowed with a sesquilinear
form $[\cdot,\cdot]$. The pair $(\V,[\cdot,\cdot])$ is
called an indefinite inner product space, and many
important concepts are associated to such a pair, some algebraic
and some topological. The combination of both is a main feature
of the general theory developed in \cite{bognar}.
The form $[\cdot,\cdot]$ defines an
orthogonality: two vectors $v,w\in \V$ are orthogonal if
$[v,w]=0$, and two linear subspaces $\V_1$ and $\V_2$ of $\V$ are
orthogonal if every vector of $\V_1$ is orthogonal to every
vector of $\V_2$. Orthogonal sums will be denoted by the symbol
$[+]$. Note that two orthogonal spaces may intersect. We will
denote by the symbol $[\oplus]$ a direct orthogonal sum. A
complex vector space $\V$ is a Krein space if it can be written
(in general in a non-unique way) as
$\V=\V_+[\oplus]\V_-$, where $(\V_+,[\cdot,\cdot])$ and $(\V_-,-[\cdot,\cdot])$ are
Hilbert spaces. When the space $\V_-$ (or, as in \cite{ikl}, the space
$\V_+$) is finite dimensional (note that this property does not depend on the decomposition), $\V$ is called a Pontryagin
space.\\

Krein spaces were introduced by  Krein and later and independently
by L. Schwartz in \cite{schwartz} where they were called
"Hermitian spaces"; for historical remarks, we refer to [4, pp.
207-209]. Besides the book of Bognar \cite{bognar}, on which is based this
work, we refer to
\cite{azih,MR1364446} for the theory of Krein spaces and of their
operators, and to \cite{ikl} for
the case of Pontryagin spaces.

Among other topics, Krein spaces appear in a natural
way in the theory of interpolation for operator-valued Schur
functions (see e.g., \cite{bh-1986,MR756761,ball-fang,MR2394102}). The
motivation for the
present work came in particular from the desire to extend
interpolation theory for operator-valued Schur functions to the case of
slice-hyperholomorphic functions,
see the forthcoming work \cite{abcs2-milan}, as we now explain.
Let $\mathcal Y$ and ${\mathcal U}$ be two Hilbert spaces. We
denote by $\mathcal S({\mathcal U},\mathcal Y)$ the class of
$\mathbf L({\mathcal U}, \mathcal Y)$-valued functions analytic
and contractive in the open unit disk. To define the
left-interpolation problem in this class we need a third Hilbert
space, say $\mathcal X$, and two operators $A\in\mathbf
L(\mathcal X)$ and $C\in\mathbf L(\mathcal X, \mathcal Y)$. We
assume that the series ${\displaystyle\sum_{n=0}^\infty A^{*n}C^*CA^n}$
converges in the strong operator topology. The interpolation
problem at hand is to find all (if any) functions $S\in\mathcal
S({\mathcal U},\mathcal Y)$ such that
\[
\sum_{k=0}^\infty A^{*k}C^*S_k=N^*,
\]
where the $S_k\in\mathbf L({\mathcal U},\mathcal Y)$ are the
coefficients of the power expansion of $S$ at the origin and $N\in\mathbf
L({\mathcal X},\mathcal U)$ is given. Krein
spaces appear as follows in the solution of this problem. Let $P$
be the solution of the Stein equation
\[
P-A^*PA=\begin{pmatrix}C^* & N^*\end{pmatrix}J\begin{pmatrix}C
\\ N\end{pmatrix}, \quad\mbox{where}\quad J=\begin{pmatrix}I_{\mathcal Y} & 0 \\ 0
&-I_{{\mathcal U}}\end{pmatrix}
\]
and assume $P$ positive and boundedly invertible. We endow the
space ${\mathcal K}=\mathcal X \oplus \mathcal Y \oplus {\mathcal
U}$ with the indefinite metric defined by
\[
\widetilde{J}=\begin{pmatrix}P & 0
\\ 0 & J\end{pmatrix}.
\]
A key result in the arguments is that the space
\[
{\mathcal K}_0: ={\rm Ran} \begin{pmatrix} A \\ C \\
N\end{pmatrix}
\]
is a closed uniformly positive subspace of $\mathcal K$ and thus it is ortho-complemented.  The proof of this last fact is  in \cite{bognar},
 and requires a long chain of preliminary results and we
are not aware of any shortcut proof. As we previously remarked, the main
purpose of this paper is to prove the counterpart of
this fact in the quaternionic setting. To this end, we first need to prove
some algebraic as well as topological results for quaternionic vector
spaces which are of independent interest.
The complex version of these results can be found in  \cite{bognar},
\cite{DS1}, \cite{Rud91} (we will give more precise references where
appropriate).
In most cases, the proofs are not
substantially different from the proofs of the corresponding
results in the complex case. However, since  we
are not aware of any reference in which these results are explicitly
proven
in the quaternionic setting, we repeat them here.\\

The paper consists of nine sections besides the introduction, and its
outline is as follows: Sections 2 and 3 are
devoted to quaternionic topological vector spaces and to some basic
functional analysis theorems in a quaternionic
setting. These sections set the framework for the following sections,
where one consider quaternionic vector
spaces endowed with a possibly degenerate and non positive inner product.
The algebraic aspects of such spaces are
studied in Sections 4, 5 and 6. The notion of fundamental decomposition is
studied in detail and plays a key role in
the paper in the later sections. Topologies which make the inner product continuous (called {\sl majorants}) are studied
in Sections 6, 7 and 8. Finally, some aspects of Krein spaces are studied
in Section 9.

\section{Quaternionic topological vector spaces}
\setcounter{equation}{0}
\label{tvs}
In this paper $\mathbb H$ denotes the algebra of real quaternions.
We send to \cite[Chapter I]{MR83k:46003} and to \cite[p. 446]{as3}
for the basic definitions of a vector space over $\mathbb H$. In
this paper we will treat the case of right quaternionic vector
spaces. The case of left quaternionic vector
spaces may be treated in an analogous way. It is also useful to
recall that if $\W$ and $\V$ are two right quaternionic vector
spaces, an operator $A\,\,:\,\, \V\longrightarrow \W$ is
right linear if
\[
A(v_1q_1+v_2q_2)=(Av_1)q_1+(Av_2)q_2,\quad\forall v_1,v_2\in
\V\quad{\rm and}\quad \forall q_1,q_2\in\mathbb H.
\]
For future reference we single out the following result, which is
true for the case of vector spaces over any field or skew field. The
claims are
\cite[Th\'eor\`eme 1 and Proposition 4, Ch. 2, \S7]{MR43:2}
respectively.

\begin{Tm}
$(a)$ Every right quaternionic vector space has a basis.\\
$(b)$ Every (right) linear subspace of a quaternionic vector space
has a direct complement. \label{tm:steinitz}
\end{Tm}

We also recall the following: If ${\mathcal V}$ is a right
quaternionic vector space and $\mathcal V_1\subset\mathcal V$ is
a (right) linear subspace of ${\mathcal V}$, the quotient space
${\mathcal V}/{\mathcal V_1}$ endowed with
\[
(v+\mathcal V_1)q=vq+\mathcal V_1
\]
is also a right quaternionic vector space. Here $v +\mathcal V_1$
denotes the equivalence class in the quotient space ${\mathcal
V}/{\mathcal V_1}$ of $v\in \mathcal V_1$.\\

Given a right quaternionic vector space $\mathcal{V}$, a
semi-norm is defined (as in the complex case) as a map
$p:\,\mathcal{V}\to \mathbb R$  such that
\begin{equation}
\label{sn1}
p(v_1+v_2)\le p(v_1)+p(v_2),\quad\forall v_1,v_2\in
\mathcal{V},
\end{equation}
and
\begin{equation}
\label{sn2} p(vc)=|c|p(v),\quad \forall v\in \mathcal{V}\,\,{\rm
and}\,\, c\in\mathbb H.
\end{equation}

\begin{Rk}{\rm Note that (\ref{sn2}) implies that $p(0)=0$ and (\ref{sn1})
implies
\[
0=p(v-v)\leq 2p(v),
\]
so that a semi-norm has values in $\mathbb{R}^+$.}
\end{Rk}

As it is well known,  see \cite[Ch. II, $\S$1]{MR83k:46003},
given a vector space over a non discrete valued division ring it
is possible to introduce the notion of semi-norm. We observe that
one can give the notion of semi-norm in the framework of modules
over a Clifford algebra, see \cite{bds}. However, in that case,
(\ref{sn2}) is required only when $c\in\mathbb{R}$ while in
general it has to be replaced by the weaker condition $p(vc)\leq
C |c| p(v)$, where $C$ is  a suitable
constant.\\

A family of semi-norms on $\mathcal V$ gives rise
to a topology which, at
least in the cases of complex or real vector spaces, leads to a locally convex space.\\

Let $p$ be a semi-norm and set
$$
{U}_{v_0}(p,\alpha)\stackrel{\rm def.}{=}\{v\in\mathcal{V}\ |\
p(v-v_0)<\alpha \}.
$$
A family $\{p_\gamma\}_{\gamma\in\Gamma}$ of semi-norms on
$\mathcal{V}$ defines a topology  on $\mathcal{V}$, in which a
subset ${U}\subseteq\mathcal{V}$ is said to be open if and only if for every
$v_0\in U$ there are $\gamma_1,\ldots
,\gamma_n\in\Gamma$ and $\varepsilon>0$ such that 
$v\in U_{v_0}(p_{\gamma_j},\varepsilon)$, $j=1,\ldots ,n$, implies $v\in  {U}$.\\

\begin{Rk}
{\rm All the spaces considered here will be right linear, and in
general we will use the terminology {\sl quaternionic vector
space} rather than {\sl right quaternionic vector space}.
Similarly we will speak of linear operators rather than right
linear operators.}
\end{Rk}

A quaternionic vector space $\mathcal{V}$ is also a vector space over
$\mathbb{R}$. It is immediate to verify
using (\ref{sn1}) and (\ref{sn2}) that when it is endowed with the
topology induced by  a family of semi-norms, it is a
locally convex space.
  \begin{Dn}
{\rm A set ${U}$ in a topological quaternionic vector space $\mathcal{V}$
is called {\em balanced} if $vc\in{U}$,
whenever $v\in{U}$ and $c\in\mathbb{H}$ with $|c|\leq 1$. A set ${U}\in\mathcal{V}$
is said to be {\em absorbing} if for any
$v\in\mathcal{V}$ there exists $c>0$ such that $vc^{-1}\in{U}$.}
  \end{Dn}
\begin{Pn}
Let $p$ be a semi-norm on a quaternionic vector space
$\mathcal{V}$, let $\alpha>0$. Then the set ${U}_0(p,\alpha)=
  \{v\in\mathcal{V}\ |\ p(v)<\alpha\}$ is balanced and absorbing.
  \end{Pn}
  \begin{proof}
  By (\ref{sn2}), if $v\in{U}_0(p,\alpha)$ and $|c|\leq 1$ then
$p(vc)=|c|p(v)<\alpha$ so ${U}_0(p,\alpha)$ is balanced.
  Similarly one proves that  ${U}_0(p,\alpha)$ is absorbing.
  \end{proof}
We recall the definition  of the  Minkowski functional $p_{{U}}$
associated to a convex, balanced and absorbing set ${U}$:
\begin{equation}
\label{Minkowski}
p_U(v)=\inf A_v\quad\mbox{where}\quad
A_v=\{a>0: \, va^{-1}\in U\}, \qquad v\in V.
\end{equation}

\begin{Pn}
  Let  $\mathcal V$ be a quaternionic vector space, and let
${U}$ be a convex, balanced, absorbing set containing $0$. Then
the Minkowski functional $p_{U}(v)$  is a semi-norm on
$\mathcal{V}$.
  \end{Pn}
  \begin{proof}
  Let $v_1,v_2\in\mathcal{V}$ and $c\in A_{v_1}$, $d\in A_{v_2}$.
  Then $v_1c^{-1}+v_2d^{-1}\in{U}$ or, equivalently, $v_1+v_2\in c\,
{U}+d\,
  {U}=(c+d){U}$  since ${U}$ is convex.
Thus $c+d\in A_{v_1+v_2}$ and $p_U(v_1+v_2)\leq c+d$, from which
we conclude \eqref{sn1} since $c$ and $d$ are arbitrary.\\
To prove (\ref{sn2}), we begin by considering $\lambda >0$ and
$v\in\mathcal{V}$.  Take any $c\in A_v$, then we have
$vc^{-1}\in{U}$ and
  $v\lambda (\lambda c)^{-1}\in{U}$ and so $\lambda c\in A_{v\lambda}$ and
then
  $p_U(v\lambda )\leq \lambda c$. By the arbitrariness of $c$ it follows
  that $p_{{U}}(v\lambda)\leq \lambda p_{{U}}(v)$. By replacing $v$ by
$v\lambda$
and $\lambda$ by $\lambda^{-1}$ we obtain $p_U(v)\leq
\lambda^{-1}p_U(v\lambda)$ or, equivalently, $\lambda p_U(v)\leq
p_U(v\lambda)$. Thus $\lambda p_U(v)= p_U(v\lambda)$.
\\
If we consider $\lambda=0$ then (\ref{sn2}) is trivial since
$p_U(0)=0$ by definition \eqref{Minkowski}.  Thus we assume
now that $\lambda\in\mathbb{H}$ and $\lambda\not=0$. Let
$v\in\mathcal{V}$ and $c\in A_v$. Since ${U}$ is balanced, then
$v c^{-1}\in{U}$ and also $v
\frac{\lambda}{|\lambda|}c^{-1}\in{U}$ and so
$p_{{U}}(v\lambda)\leq |\lambda| c$. Since $c$ is arbitrary, we
have $p_{{U}}(v\lambda)\leq |\lambda| p_U(v)$. The reverse
inequality is obtained by  replacing $v$ by $v\lambda$ and
$\lambda$ by $\lambda^{-1}$. The statement follows.
\end{proof}
 \begin{Pn}
 A topological quaternionic vector space is locally convex if and only
if the topology is defined by a family of semi-norms.
 \end{Pn}
\begin{proof}
The "if" part of the statement has already been discussed. To show the
"only if" part, consider a base $B$ of neighborhood at $0$ consisting of convex and balanced
open sets. Since the multiplication by a scalar on right is continuous, each $U\in B$ is absorbing. Then for $U\in B$ we define
$p_U(v)=\inf  A_v$ (see (\ref{Minkowski})) and so $p_U$ is the Minkowski
functional. The family $\{p_{{U}}\}_{U\in B}$ is then a family of semi-norms  such that
$$
\{v\in\mathcal{V}\ : \ p_U(v)<1\}\subseteq U \subseteq \{v\in\mathcal{V}\ : \ p_U(v)\leq 1\}
$$
and the statement follows.
\end{proof}

We conclude this section by mentioning that the topology induced by the family of
semi-norms $\{p_\gamma\}_{\gamma\in\Gamma}$ is
 Hausdorff  if and only if the condition
$p_\gamma(v)=0$ for all $\gamma\in\Gamma$ implies $v=0$.

\section{Principles of quaternionic functional analysis }
\setcounter{equation}{0}
The material in this section is classical for complex Fr\'echet spaces and
can be found e.g. in \cite[Chapter II]{DS1} or \cite[Chapter 2]{Rud91}.\\
Let $\mathcal{V}$ be a quaternionic Fr\'echet space, that is a quaternionic
locally convex topological vector space which is metrizable and
complete, and let $\rho$ be an associated metric.
 For the sake of simplicity, in the sequel we will write
 $|u-w|$ instead of $\rho(u,w)$.\\

 We now prove a result for continuous (not necessarily linear) maps
 which implies  the principle of uniform boundedness.
\begin{Tm}
For each $a\in A$, where $A$ is a set, let $S_a$ be a continuous
map of a quaternionic Fr\'echet space $\mathcal{V}$ into a
quaternionic Fr\'echet space $\mathcal{W}$, which satisfies the
following properties
\begin{itemize}
\item[(a)] $|S_a(u+w)|\leq |S_a(u)|+|S_a(w)|$,  $\forall u,w\in
\mathcal{V}$,
\item[(b)] $|S_a(w\alpha)|=|S_a(w)\alpha|$,  $\forall w\in  \mathcal{V}$,
$\forall \alpha \geq 0$.
\end{itemize}
If, for each $u\in \mathcal{V}$, the set $\{ S_av\}_{a\in A}$ is
bounded, then $\lim_{v\to 0}S_av=0$ uniformly in $a\in A$.
\label{T:1}
\end{Tm}
\begin{proof}
For $\varepsilon>0$, $a\in A$ and a positive integer $k$, the set
$$
 \mathcal{V}_k\stackrel{\rm def.}{=}\left\{u\in  \mathcal{V}\ : \
\left|\frac{1}{k}S_a(u)\right|+\left|\frac{1}{k}S_a(-u)\right|\leq
\frac{\varepsilon}{2}\right\}
$$
is closed since $S_a$ are continuous.
Moreover, by assumption, the sets $\{ S_av\}_{a\in A}$ are  bounded, so
$$
\mathcal{V}=\bigcup_{k=1}^\infty\mathcal{V}_k.
$$
By the Baire category theorem, there
exists a $\mathcal{V}_{k_0}$ that contains a ball $B(v_0,\delta)$
with center at $v_0$ and radius $\delta >0$. Let $|u|<\delta$. Then
$v_0$ and $v_0+u$ both belong to $B(v_0,\delta)$
and so they both are in $\mathcal{V}_{k_0}$. Thus we
have
$$
\left|\frac{1}{k_0}S_a(v_0+u)\right|\leq
\frac{\varepsilon}{2}\qquad {\rm and}\qquad
 \left|\frac{1}{k_0}S_a(-v_0)\right|\leq \frac{\varepsilon}{2}.
$$
Using assumption (a) we deduce
$$
\left|\frac{1}{k_0}S_a(u)\right|\leq \left|\frac{1}{k_0}
S_a(v_0+u)\right|+\left|\frac{1}{k_0}S_a(-v_0)\right|,
$$
and using assumption (b) we get
$$
\left|\frac{1}{k_0}S_a(u)\right|=\Big|S_a\Big(\frac{1}{k_0}u\Big) \Big|
\leq \varepsilon,\ \ \ |u|<\delta,\ \ a\in A.
$$
Now observe that the mapping $v\mapsto v/k_0$ is a homeomorphism
of $ \mathcal{V}$ into itself since $\mathcal{V}$ is a
topological vector space and thus the multiplication by a
quaternionic scalar (in particular, real) is continuous. Thus
$\lim_{v\to 0}S_av=0$ uniformly in $a\in A$.
\end{proof}
In the case of linear maps, Theorem \ref{T:1} gives the following result which
 will be used in the proofs of
Theorem \ref{Mairie-de-Montreuil-ligne-9} and Proposition \ref{pn6.2} below.

\begin{Tm}[Principle of uniform boundedness]
For each $a\in A$, where $A$ is a set, let $T_a$ be continuous
linear map of a quaternionic Fr\'echet space $\mathcal{V}$ into a
quaternionic Fr\'echet space $\mathcal{W}$. If, for each $u\in
\mathcal{V}$, the set $\{ T_av\}_{a\in A}$ is bounded, then
$\lim_{v\to 0}T_av=0$ uniformly in $a\in A$. \label{tmunivb}
\end{Tm}
For further reference, we repeat also the formulation of the principle
of uniform boundedness for quaternionic Banach spaces.
\begin{Tm}
\label{PrinUniBoun} Let $\mathcal{V}$ and $\mathcal{W}$ be two
quaternionic Banach spaces and let $\{T_a\}_{a \in A}$ be bounded
linear maps from  $\mathcal{V}$ to $\mathcal{W}$. Suppose that
${\displaystyle\sup_{a \in A}\|T_\alpha v\| <\infty}$ for any
$v\in\mathcal{V}$.
Then
$$
\sup_{a \in A}\|T_\alpha\| <\infty.
$$
\end{Tm}
The next result is the
quaternionic counterpart of the open mapping theorem:

\begin{Tm}[Open mapping theorem]\label{OPNMAP}
Let $\mathcal{V}$ and $\mathcal{W}$ be two quaternionic Fr\'echet
spaces, and let $T$ be a  linear continuous quaternionic map
from $\mathcal{V}$ onto $\mathcal{W}$. Then the image of every
open set is open.
\end{Tm}

\begin{proof} Let $B_{\mathcal{V}}(r)\subset\mathcal{V}$ denote
the open ball of radius $r>0$ and centered at the origin and let
$B_{\mathcal{V}}(r)- B_{\mathcal{V}}(r)$ be the set of elements of the form
$u-v$ where $u,v\in B_{\mathcal{V}}(r)$. Since  the function $u-v$ is continuous in  $u$
and $v$, there exists a ball $B_{\mathcal{V}}(r')$, for suitable
$r'>0$, such that $B_{\mathcal{V}}(r')-
B_{\mathcal{V}}(r')\subseteq B_{\mathcal{V}}(r)$. For every $v\in \mathcal{V}$ we have
that $v/n\to 0$ as $n\to
\infty$ so $v\in nB_{\mathcal{V}}(r')$ for a suitable $n\in
\mathbb{N}$. So
$$
\mathcal{V}=\bigcup_{n=1}^\infty n B_{\mathcal{V}}(r')\quad
{\rm and}\quad\mathcal{W}=T\mathcal{V}=\bigcup_{n=1}^\infty n
TB_{\mathcal{V}}(r').
$$
By the Baire category theorem one of the closures $\overline{n
TB_{\mathcal{V}}(r')}$ contains a non empty open set. The map
$w\mapsto nw$ is a homeomorphism in $\mathcal{W}$ and
$\overline{TB_{\mathcal{V}}(r')}$ contains a non empty open set
denoted by $\mathcal{B}$, so
\begin{equation}
\overline{TB_{\mathcal{V}}(r)}\supseteq
\overline{TB_{\mathcal{V}}(r')- TB_{\mathcal{V}}(r')}\supseteq
\overline{TB_{\mathcal{V}}(r')}- \overline{TB_{\mathcal{V}}(r')}
\supseteq \mathcal{B}- \mathcal{B}.
\label{3}
\end{equation}
The map $w\mapsto u-w$ is a homeomorphism and hence  the
set $u- B_{\mathcal{V}}(r)$ is open. Since the set $\mathcal{B}-
\mathcal{B}=\bigcup_{u\in \mathcal{B}}(u-\mathcal{B})$ is open
(as the union of open sets)  and contains the origin, we conclude from
\eqref{3} that $\overline{TB_{\mathcal{V}}(r)}$ contains a
neighborhood of the origin.

\smallskip

Fix an arbitrary $\varepsilon_0$ and let
$\varepsilon_\ell >0$ be a sequence such that $\sum_{\ell\in
\mathbb{N}}  \varepsilon_\ell <\varepsilon_0$. Then there exists
a sequence $\theta_\ell >0$ with $\theta_\ell\to 0$ such that
\begin{equation}\label{overTBV}
\overline{TB_{\mathcal{V}}(\varepsilon_\ell)}\supset
B_{\mathcal{W}}(\theta_\ell), \ \ \ \ell\in \mathbb{N}\cup\{0\}.
\end{equation}
We now take an arbitrary $w\in B_{\mathcal{W}}(\theta_0)$ and show that
$w=Tv$ for some $v\in B_{\mathcal{V}}(2\varepsilon_0)$.
To this purpose we follow a recursive procedure. From (\ref{overTBV}) for $\ell=0$
there exists $v_0\in  B_{\mathcal{V}}(\varepsilon_0)$ such that $|w-Tv_0|<\theta_1$.
Since $w-Tv_0\in B_{\mathcal{W}}(\theta_1)$, then $w-Tv_0\in TB_{\mathcal{V}}(\varepsilon_1)$  and  again from
(\ref{overTBV}) with $\ell=1$, there exists $v_1\in
B_{\mathcal{V}}(\varepsilon_{1})$ such that $|w-Tv_0-Tv_1 |<\theta_2$.
Iterating this procedure, we construct a sequence $\{v_n\}_{n\in \mathbb{N}}$ such that
$v_n\in B_{\mathcal{V}}(\varepsilon_n)$ and
\begin{equation}\label{differencevTxi}
|w-T\sum_{\ell=0}^nv_\ell  |<\theta_{n+1},\ \ \ n\in
\mathbb{N}\cup\{0\}.
\end{equation}
Let us denote $p_m=\sum_{\ell=0}^mv_\ell$. Then $\{p_m\}$
is a Cauchy sequence since
$$
|p_m-p_n|=|v_{n+1}+...+v_m|<\varepsilon_{n+1}+\ldots +\varepsilon_m\quad\mbox{for}\quad
m>n.
$$
Therefore the series
$\sum_{\ell=0}^\infty v_\ell $ converges to  a point $v\in\mathcal V$ with
$|v|\leq \sum_{\ell=0}^\infty \varepsilon_{\ell}=2 \varepsilon_0$.
Since $T$ is continuous, we conclude from  (\ref{differencevTxi})
that $w=Tv$. We thus showed that an arbitrary ball
$B_{\mathcal{V}}(2\varepsilon_0)$ is mapped onto the set
$TB_{\mathcal{V}}(2\varepsilon_0)$ which contains the ball
$B_{\mathcal{W}}(\theta_0)$. So if $\mathcal{N}$ is  a
neighborhood of the origin in $\mathcal{V}$ then $T\mathcal{N}$
contains a neighborhood of the origin of $\mathcal{W}$. Since $T$
is linear then the above procedure works for every neighborhood
of every point.
\end{proof}
\begin{Tm}[Banach continuous inverse  theorem]\label{BCInverse}
Let $\mathcal{V}$ and $\mathcal{W}$ be two quaternionic Fr\'echet
spaces and let $T: \; \mathcal{V}\to \mathcal{W}$ be a one-to-one linear continuous
quaternionic map.  Then $T$ has a  linear continuous inverse.
\end{Tm}
\begin{proof}
By Theorem \ref{OPNMAP} $T$ maps open sets
onto open sets, so if we write $T$ as $(T^{-1})^{-1}$, it is
immediate that $T^{-1}$ is continuous. Now take $w_1$, $w_2\in
\mathcal{W}$ and $v_1$, $v_2\in \mathcal{V}$ such that
$Tv_1=w_1$, $Tv_2=w_2$ and $p\in \mathbb{H}$. Then
$$
T(v_1+v_1)=Tv_1+Tv_2=w_1+w_2, \ \ \ T(v_1p)=T(v_1)p=w_1p
$$
and hence
$$
T^{-1}(w_1+w_2)=v_1+v_2\quad\mbox{and}\quad
T^{-1}(w_1p)=v_1p,
$$
so $T^{-1}$ is  linear quaternionic operator.
\end{proof}
\begin{Dn}
Let $\mathcal{V}$ and $\mathcal{W}$ be two quaternionic Fr\'echet
spaces. Suppose that $T$ is a  quaternionic operator whose
domain
 $\mathcal{D}(T)$ is a linear manifold contained in $\mathcal{V}$ and
whose
 range belongs to $\mathcal{W}$.
The graph of $T$ consists of all point $(v,Tv)$, with  $v\in
\mathcal{D}(T)$, in the product space $\mathcal{V}\times
\mathcal{W}$.
\end{Dn}
\begin{Dn}
We say that $T$ is a closed operator if its graph is closed in
$\mathcal{V}\times \mathcal{W}$.
\end{Dn}
\begin{Rk}{\rm
Equivalently we can say that $T$ is closed if $v_n\in
\mathcal{D}(T)$, $v_n\to v$, $Tv_n\to y$ imply that $v\in
\mathcal{D}(T)$ and $Tv=y$.}
\end{Rk}

The following theorem can be found also in \cite[Corollaire 5, p.  I.19]{MR83k:46003}.

\begin{Tm}[Closed graph theorem]\label{closedGR}
Let $\mathcal{V}$ and $\mathcal{W}$ be two quaternionic Fr\'echet
spaces. Let $T: \mathcal{V}\to \mathcal{W}$ be a linear closed
quaternionic operator. Then $T$ is continuous.
\end{Tm}

\begin{proof}
Since $\mathcal{V}$ and $\mathcal{W}$ are two quaternionic
Fr\'echet spaces
 we have that $\mathcal{V}\times \mathcal{W}$ with the distance
$|(v,w)|_{\mathcal{V}\times
\mathcal{W}}=|v|_{\mathcal{V}}+|w|_{\mathcal{W}}$  is a  quaternionic
Fr\'echet space.
The graph of $T$ denoted by $\mathcal{G}(T)=\{(v,Tv),\  v\in
\mathcal{D}(T)\}$ is a closed linear manifold in the product space
$\mathcal{V}\times \mathcal{W}$ so it is a quaternionic Fr\'echet
space. The projection
$$
P_{\mathcal{V}}:\mathcal{G}(T)\mapsto \mathcal{V},\ \ \ \
P_{\mathcal{V}}(v,Tv)= v
$$
is one-to-one and onto, linear and continuous so by Theorem
\ref{BCInverse} its inverse $P^{-1}_{\mathcal{V}}$ is continuous.
Now consider the projection
$$
P_{\mathcal{W}}:\mathcal{G}(T)\mapsto \mathcal{W},\ \ \ \
P_{\mathcal{W}}(v,Tv)= Tv,
$$
observing that $T=P_{\mathcal{W}}P^{-1}_{\mathcal{V}}$ we get the
statement.
\end{proof}

\section{Ortho-complemented spaces}
\setcounter{equation}{0}
From this section on, we focus on quaternionic vector spaces endowed
with an inner product, defined as follows:
\begin{Dn}\label{innerproduct}
Let ${\mathcal V}$ be a quaternionic vector space. The map
\[
[\cdot, \cdot]\,\, :\,\, {\mathcal V}\times {\mathcal
V}\quad\longrightarrow\quad \mathbb H
\]
is called an inner product if it is a (right) sesquilinear form:
\[
[v_1c_1,v_2c_2]=\overline{c_2}[v_1,v_2]c_1,\quad \forall
v_1,v_2\in {\mathcal V},\,\,{and}\,\, c_1,c_2\in\mathbb H,
\]
and Hermitian:
\[
[v,w]=\overline{[w,v]},\quad \forall v,w\in\mathbb {\mathcal V}.
\]
\label{dn:hermi}
\end{Dn}

We will call the pair $({\mathcal V},[\cdot,\cdot])$ (or the
space ${\mathcal V}$ for short when the form is understood from
the context) a (right) quaternionic indefinite inner product
space. A form is called positive (or non-negative) if
$[v,v]\ge 0$ for all $v\in {\mathcal V}$.
\begin{Rk}{\rm
Note that the Cauchy-Schwarz inequality holds for positive inner
product spaces; see  \cite[Lemma 2.2, p. 5]{bognar} for the
classical case and \cite[Lemma 5.6 and Remark 5.7, p. 447]{as3}
and the references therein for the quaternionic case. Multiplying
the inner product by $-1$, we see that the Cauchy-Schwarz
inequality holds in inner product spaces for which the inner
product is negative. We will call an inner product space definite, if it is either positive or negative.}
\end{Rk}

%
%
The definitions on indefinite product spaces over $\mathbb
C$ reviewed in the introduction carry over when one considers the
quaternions. In particular,
two elements $v$ and $w$  in ${\mathcal V}$ will be called
orthogonal if $[v,w]=0$, and two vector subspaces ${\mathcal
V}_1$ and ${\mathcal V}_2$ of ${\mathcal V}$ are orthogonal if
every element of ${\mathcal V}_1$ is orthogonal to every element
of ${\mathcal V}_2$. Two orthogonal subspaces ${\mathcal V}_1$
and ${\mathcal V}_2$ may have a non trivial intersection. When
their intersection reduces to the zero vector we denote  by
${\mathcal V}_1{[\oplus]}{\mathcal V}_2$ their direct orthogonal sum.\\

For $\mathcal L\subset {\mathcal V}$ we set
\[
\mathcal L^{[\perp]}=\left\{v\in {\mathcal V}\,\,:\,\
[v,w]=0,\,\,\forall w\in\mathcal L\right\}.
\]
The definition makes sense even when $\mathcal L$ is not a linear
space, but a mere subset of ${\mathcal V}$, and the set $\mathcal
L^{[\perp]}$ is always a linear space. It is called the {\sl
orthogonal companion} of $\mathcal L$. Note that
\begin{equation}
\label{atlanta-ineql}
\mathcal L\subset\left(\mathcal
L^{[\perp]}\right)^{[\perp]}\stackrel{\rm def.}{=} \mathcal
L^{[\perp\perp]}.
\end{equation}

A linear subspace $\mathcal L$  is called non-degenerate if its {\sl
isotropic part} $\mathcal L^0\stackrel{\rm def.}{=}\mathcal
L\cap\mathcal L^{[\perp]}$ is trivial.

\begin{Pn}
Let $\V$ be a quaternionic inner product space, and let $\V^0$ be
its isotropic part. The formula
\begin{equation}
\label{innrerperp} [v+\mathcal V^0, w+\mathcal
V^0]_q\stackrel{\rm def.}{=}[v,w]
\end{equation}
defines a non-degenerate indefinite inner product on $\mathcal
V/\mathcal V^0$.
\end{Pn}

\begin{proof}
It suffices to note that formula \eqref{innrerperp} is well
defined (that is, does not depend on the specific choice of $v$
and $w$).
\end{proof}


We now move to ortho-complemented spaces.
Following \cite[p. 18]{bognar} we say that the space $\mathcal L$
is {\sl ortho-complemented} if $\mathcal V$ is spanned by
$\mathcal L$ and $\mathcal L^{[\perp]}$. As explained in the
introduction, the motivation for the present work was to prove in
the quaternionic setting that a uniformly positive and closed 
subspace of a quaternionic Krein space is ortho-complemented (the
definitions of these various notions appear in the sequel and the
result itself is Theorem \ref{tm:sofsof} below). In this section
we will prove some results in the quaternionic setting whose
counterparts for complex vector spaces can be found in
\cite[Chapter I]{bognar}. We begin by stating the following direct
consequence
of \eqref{atlanta-ineql}.

\begin{Pn}
\label{pn-ortho} Let ${\mathcal V}$ be a quaternionic vector
space, and let $\mathcal L$ denote a linear  subspace of
$\mathcal V$ which is ortho-complemented. Then, $\mathcal
L^{[\perp]}$ is also ortho-complemented.
\end{Pn}

A linear subspace ${\mathcal M}\subset {\mathcal V}$ is called
positive if $[m,m]\ge 0$ for all $m\in {\mathcal M}$. It is
called strictly positive if the inequality is strict for all
$m\not=0$. Similar definitions hold for negative and strictly
negative subspaces. A linear subspace will be called definite if
it is either positive or negative, and indefinite otherwise. It
will be called neutral if $[m,m]=0$ for all $m\in {\mathcal M}$.\\


The following
result is stated for future reference. Note that in the statement
the spaces may have a non trivial intersection.

\begin{Pn}\label{Prop52}
Let $({\mathcal V},[\cdot,\cdot])$ be an indefinite quaternionic
vector
space.\\
$(a)$ Let ${\mathcal M}_1,\ldots, {\mathcal M}_n$ be $n$ pairwise
orthogonal subspaces of ${\mathcal V}$. Assume that all the
${\mathcal M}_i$ are positive (resp. neutral, negative, strictly
positive, strictly negative). Then the space spanned by
the ${\mathcal M}_i$ has the same property.\\
$(b)$ Let $m_1,\ldots , m_n\in {\mathcal V}$ be vectors which are
positive (resp. neutral, negative, strictly positive, strictly
negative). Then, for every choice of $q_1,\ldots, q_n\in\mathbb
H$ the vector
\[
m=\sum_{j=1}^n m_jq_j
\]
is positive (resp. neutral, negative and when at least one of the
$q_j\not=0$, strictly positive, strictly negative).
\end{Pn}

\begin{proof}\mbox{}\\
$(a)$ An element $m$ is in the linear span of ${\mathcal
M}_1,\ldots,{\mathcal M}_n$ if and only if it can be written (in
general in a non-unique way) as
\begin{equation}
\label{eq-m}
m=\sum_{j=1}^n m_j q_j,
\end{equation}
where the $m_j\in {\mathcal M}_j$ for $j=1,\ldots n$. Then,
\[
\begin{split}
[m,m]&=[\sum_{j=1}^n m_j q_j,\sum_{k=1}^n m_kq_k]\\
&=\sum_{j=1}^n [m_j,m_j] |q_j|^2+\sum_{j\not
=k}\bar q_k[m_j,m_k]q_j\\
&=\sum_{j=1}^n [m_j,m_j]|q_j|^2
\end{split}
\]
since by hypothesis, $[m_j,m_k]=0$ for $j\not=k$. The result
follows.\\

$(b)$ This item follows from the fact that the linear span of a positive vector
is a one dimensional subspace which is positive, and similarly
for the other cases at hand.
\end{proof}

We now  briefly discuss some properties of the isotropic part of
an indefinite quaternionic inner product space.

\begin{Pn} 
Let ${\mathcal V}$ denote a definite quaternionic inner product
space.
Then:\\
$(a)$ Assume ${\mathcal V}$ positive (resp. negative). Then, an
element $v$ belongs to the isotropic part ${\mathcal V}^0$ of
${\mathcal V}$ if
and only if it is neutral: $[v,v]=0$.\\
$(b)$ Assume ${\mathcal V}$ neutral. Then, the inner product
vanishes identically in ${\mathcal V}$.\\
$(c)$ A neutral subspace of $\mathcal N\subset \mathcal V$ is
ortho-complemented if and only if it is included in the isotropic
part of $\mathcal V$.
\label{cs1234}
\end{Pn}

\begin{proof} The first two statements are direct consequences of the
Cauchy-Schwarz inequality, which, as already remarked, holds in
definite quaternionic inner product spaces. As for the third
claim, item $(b)$ implies that $\mathcal N\subset\mathcal
N^{[\perp]}$. Thus, $\mathcal N$ is ortho-complemented if and
only if $\mathcal V=\mathcal N^{[\perp]}$, that is $\mathcal N$
is orthogonal to $\mathcal V$, which is the claimed inclusion.
\end{proof}

\begin{Pn} 
Let ${\mathcal V}$ denote a quaternionic inner product space, and
let ${\mathcal M}_1,\ldots, {\mathcal M}_n$ be subspaces of
${\mathcal V}$ which pairwise are orthogonal and have
intersection reducing to $\left\{0\right\}$. Then,
\begin{equation}
\label{sumpn42} \left({[\oplus]}_{j=1}^n {\mathcal
M}_j\right)^0={[\oplus]}_{j=0}^n{\mathcal M}_j^0,
\end{equation}
where we recall that the symbol ${[\oplus]}$ denotes direct and
orthogonal sum. \label{pn4-2}
\end{Pn}

Note that, since ${\mathcal M}_j^0\subset {\mathcal M}_j$, the
sum on the right side of
\eqref{sumpn42} is indeed both direct and orthogonal.\\

\begin{proof}[Proof of Proposition \ref{pn4-2}]
Let $m$ and $\ell$ be in the (direct and orthogonal) sum of the
${\mathcal M}_j$. They can be written (in a unique way) in the
form \eqref{eq-m}:
\[
m=\sum_{j=1}^n m_j\quad{\rm and}\quad \ell =\sum_{j=1}^n \ell_j,
\]
where $m_j$ and $\ell_j$ belong to ${\mathcal M}_j$, $j=1,2,\ldots,
n$. Thus
\[
[m,\ell]=\sum_{j=1}^n[m_j,\ell_j].
\]
Thus, $m$ is orthogonal to all elements in ${\mathcal M}$ if and
only if
\[
[m_j,\ell_j]=0,\quad\forall \ell_j\in {\mathcal M}_j,\quad
j=1,2,\ldots, n,
\]
that is, if and only if $\ell$ belongs to
${[\oplus]}_{j=0}^n{\mathcal M}_j^0$.
\end{proof}

In the statement of the following proposition, the existence of a
direct complement is insured by Theorem \ref{tm:steinitz}.

\begin{Pn}
Let ${\mathcal V}$ denote a quaternionic inner product space, and let
${\mathcal V}^0$ be its isotropic part. Let ${\mathcal V}_1$ be a
direct complement of ${\mathcal V}^0$. Then ${\mathcal V}_1$ is
non-degenerate and we have the direct sum decomposition
\begin{equation}
\label{ksu-nov-2012} {\mathcal V}={\mathcal V}^0
{[\oplus]}{\mathcal V}_1.
\end{equation}
\label{lemma:5.1}
\end{Pn}

\begin{proof} Let $v\in{\mathcal V}^0\cap {\mathcal V}_1$ such that
\[
[v,v_1]=0,\quad\forall v_1\in {\mathcal V}_1.
\]
On the other hand, by definition of the isotropic part,
\[
[v,v_0]=0,\quad\forall v_0\in {\mathcal V}^0.
\]
Since ${\mathcal
V}_1$ is a direct complement of  ${\mathcal V}^0$ in $\mathcal V$,  we have $v\in {\mathcal V}^0$, and so $v=0$
since ${\mathcal V}^0\cap {\mathcal V}_1=\left\{0\right\}$. The
equality \eqref{ksu-nov-2012} follows.
\end{proof}

We now gather in form of a proposition \cite[Lemmas 6.2, 6.3 and
6.4 p. 13]{bognar}. As is remarked in \cite[p. 13]{bognar}, the
claims $(b)$ and $(c)$ in the proposition are not consequences one
of the other. Note that on page 13 of that reference, the
proof of Lemma 6.3 is in fact the proof of Lemma 6.4. Note also
that we get three other claims when replacing positive by
negative in the statements.

\begin{Pn}
\mbox{}\\
$(a)$ Let ${\mathcal V}={\mathcal V}_1[\oplus]{\mathcal V}_2$
denote an orthogonal direct decomposition of the indefinite inner
product quaternionic vector space ${\mathcal V}$, where
${\mathcal V}_1$ is positive and ${\mathcal V}_2$ is maximal
strictly negative. Then, ${\mathcal V}_1$ is maximal
positive.\\
$(b)$ The space orthogonal to a maximal positive subspace is
negative.\\
$(c)$ The space orthogonal to a maximal strictly positive subspace
is negative.
\label{lemma:6.4}
\end{Pn}

\begin{proof}\mbox{}\\
$(a)$ Let $\mathcal W_1\supset {\mathcal V}_1$ be a positive
subspace of ${\mathcal V}$ containing ${\mathcal V}_1$, let $v\in
\mathcal W_1\setminus {\mathcal V}_1$, and write $v=v_1+v_2$,
where $v_1\in {\mathcal V}_1$ and $v_2\in {\mathcal V}_2$. Then,
$v_2=v-v_1\in \mathcal W_1$ since $\mathcal W_1$ is a subspace. On
the other hand, $v_2\not =0$ (otherwise $v\in {\mathcal V}_1$)
and so $[v_2,v_2]<0$. This
contradicts the fact that $\mathcal W_1$ is positive.\\

$(b)$ Let $\mathcal L$ be a maximal positive subspace of
${\mathcal V}$, and
let $v\in \mathcal L^{[\perp]}$. We distinguish three cases:\\
\begin{enumerate}

\item If $v\not\in\mathcal L$  and $[v,v]=0$, there is nothing to
prove.

\item If $v\not\in\mathcal L$  and $[v,v]>0$,
then the space spanned by $v$ and $\mathcal L$ is positive,
contradicting the maximality of $\mathcal L$. So $[v,v]\le 0$.

\item If $v\in\mathcal L$. Then, $v\in\mathcal L\cap \mathcal
L^{[\perp]}$, and so $[v,v]=0$, which is what we wanted to prove.\\

\end{enumerate}

$(c)$ Let now $\mathcal L$ be a maximal positive definite subspace
of ${\mathcal V}$, and let $v\in \mathcal L^{[\perp]}$, different
from $0$. If $[v,v]\le 0$ there is nothing to prove. If
$[v,v]>0$, the space spanned by $v$ and $\mathcal L$ is strictly
positive, contradicting the maximality of $\mathcal L$.
\end{proof}

\begin{Pn}
Let ${\mathcal V}$ denote a quaternionic inner product space, and
let $\mathcal L={[\oplus]}_{j=1}^N \mathcal L_j$ be the direct
orthogonal sum of $\mathcal L_1,\ldots, \mathcal L_N$. Then,
$\mathcal L$ is ortho-complemented if and only if each of the
$\mathcal L_j$ is ortho-complemented.
\end{Pn}

\begin{proof}
Assume first that $\mathcal L$ is ortho-complemented, and
let $\mathcal M\subset\mathcal L^{[\perp]}$ be such that
\[
\mathcal V=\mathcal L[+]\mathcal M,
\]
where the sum is orthogonal, but need not be direct. Thus
\[
\mathcal V=({[\oplus]}_{j=1}^N\mathcal L_j)+\mathcal M.
\]
For a given $j\in\left\{1,\ldots, N\right\}$, the space
\[
\mathcal M_j=({[\oplus]}_{\substack{k=1\\k\not=j}}^N \mathcal
L_k)+\mathcal M
\]
is inside $\mathcal L_j^{[\perp]}$ and such that
$\mathcal V=\mathcal L_j[+]\mathcal M_j$. Thus, $\mathcal L_j$ is
ortho-complemented.\\

Conversely (and here we follow the proof of \cite[Theorem 8.5, p. 17]{bognar}),
assume that $\mathcal L$ is
ortho-complemented. Let $v\in\mathcal V$. For $j=1,\ldots, N$ we have
\begin{equation}
\label{decompv}
v=\ell_j+m_j,\quad {\rm with}\quad\ell_j\in\mathcal L_j\quad{\rm and}\quad
m_j\in\mathcal L_j^{[\perp]}.
\end{equation}
Let $\ell=\sum_{j=1}^N\ell_j\in\mathcal L$, and let, for $j=1,\ldots, N$
\[
w_j=m_j-\left(\sum_{\substack{k=1\\ k\not=j}}^N\ell_k\right).
\]
Let $j_1,j_2\in\left\{1,\ldots, N\right\}$,
\[
\begin{split}
w_{j_1}-w_{j_2}&=m_{j_1}-m_{j_2}-\left(\sum_{\substack{k=1\\ k\not=j_1}}^N\ell_k\right)+
\left(\sum_{\substack{k=1\\ k\not=j_2}}^N\ell_k\right)\\
&=m_{j_1}-m_{j_2}-m_{j_2}+m_{j_1}\\
&=0,
\end{split}
\]
in view of \eqref{decompv}. Thus, $w_j$ is independent of $j$. We
set $w_j=w$. We have $v=\ell+w$. Furthermore, by its very
definition, it is orthogonal to every $\mathcal L_j$, and hence
orthogonal to $\mathcal L$, and this concludes the proof.
\end{proof}

If $\mathcal L$, $\mathcal V_1$ and $\mathcal V_2$ are subspaces
of the  quaternionic vector space $\mathcal V$, and
$\mathcal V_1\subset\mathcal V_2$, then one can define a map $I$
from $\mathcal L/\mathcal V_1$ into $\mathcal V/\mathcal V_2$ via
\begin{equation}
\label{mapI}
I(\ell+\mathcal V_1)=\ell+\mathcal V_2,
\end{equation}
since $\ell\in \mathcal V_1$ implies that $\ell\in\mathcal V_2$.
In general the map $I$ will not be one-to-one. If $\ell+\mathcal
V_2=\ell^\prime+\mathcal V_2$ where $\ell$ and $\ell^\prime$
belong to $\mathcal L$, then $\ell-\ell^\prime\in \mathcal
V_2\cap\mathcal L$. This need not imply that
$\ell-\ell^\prime\in\mathcal V_1$ since  we do not have in general
\begin{equation}
\mathcal V_2\cap\mathcal L\subset\mathcal V_1.
\label{131112b}
\end{equation}

We also note the following:
\begin{equation}
I(\mathcal L^{[\perp]})=\left(I(\mathcal L)\right)^{[\perp]}
\label{ortho}
\end{equation}
where we denote by the same symbol orthogonality with respect to
the original inner product and with respect the inner product
\eqref{innrerperp}.\\

We now prove the counterpart of \cite[Theorem 9.4, p.
18]{bognar}.

\begin{Tm}
Let ${\mathcal V}$ denote a quaternionic inner product space.
Then the subspace $\mathcal L$ is ortho-complemented if and only
if the following two conditions are in force:\\
$(a)$  The isotropic part of $\mathcal L$ is included in the
isotropic part of $\mathcal V$.\\
$(b)$ The image under the map $I$ (defined by \eqref{mapI}) of
the quotient space $\mathcal L/\mathcal L^0$ is
ortho-complemented in $\mathcal V/\mathcal V^0$.
\label{tm-kansas}
\end{Tm}

\begin{proof}
We first assume that $\mathcal L$ is ortho-complemented, that is $\mathcal V=\mathcal L[+]\mathcal L^{[\perp]}$.

The inner product \eqref{innrerperp} preserves orthogonality, and
thus
\[
\mathcal V/\mathcal V_0=(\mathcal L/\mathcal V_0)[+](\mathcal
L^{[\perp]}/\mathcal V_0).
\]
We now show that the map $I$ is one-to-one and so
\[
(\mathcal L/\mathcal V_0)=I(\mathcal L/\mathcal L_0),
\]
and this will conclude the proof of the direct assertion. Every
$v$ in $\mathcal V$ can be written as
\[
v=\ell+m,\quad \ell\in\mathcal L,\quad m\in\mathcal L^{[\perp]}.
\]
Let now $\ell_0\in\mathcal L^0$. We have
\[
[\ell_0,v]=[\ell_0,\ell]+[\ell_0,m]=0
\]
and thus $\mathcal L^0\subset\mathcal V^0$.  Equation
\eqref{131112b} becomes
\begin{equation}
\label{131112c}
\mathcal V^0\cap\mathcal L\subset\mathcal L^0,
\end{equation}
which always holds, and by the discussion before the theorem the
map $I$ well defined and one-to-one and so $(b)$ holds.\\

Conversely we assume now that $(a)$ and $(b)$ hold. We prove
that $\mathcal L$ is ortho-complemented. Taking $\mathcal
V_1=\mathcal L^0$ and $\mathcal V_2=\mathcal V^0$, $(a)$ insures
that the map $I$ is well defined and equation \eqref{131112c}
holds by definition of $\mathcal V_0$. Thus the map $I$ is
one-to-one. Using $(b)$ we see that for every $v\in\mathcal V$
there exist $\ell\in\LL$ and $m\in\LL^{[\perp]}$ such that
\[
v+\mathcal V^0=\ell+\V^0+m+\V^0.
\]
Thus we have $v=\ell+m+v_0$. This concludes the proof since
$\V^0\subset\LL^0\subset\LL^{[\perp]}$.
\end{proof}

We conclude this section with results pertaining to a
non-degenerate space (that is, when $\V^0$ is trivial), and which
are corollaries of the previous discussion.

\begin{Pn}\label{corollary9.5}
Let $\V$ be a  quaternionic non-degenerate inner product
space. Then:\\
$(a)$ Every ortho-complemented subspace is non-degenerate.\\
$(b)$ Let $\LL\subset \V$ be ortho-complemented. Then
$\LL={\LL^{[\perp \perp]}}$.
\end{Pn}

\begin{proof}\mbox{}
$(a)$ follows directly from  Theorem
\ref{tm-kansas} $(a)$ since
\[
\LL\cap\LL^{[\perp]}\subset\V^0=\left\{0\right\}.
\]
As for $(b)$, we always have
\begin{equation}
\label{161112ineq} \LL\subset{\LL^{[\perp \perp]}}.
\end{equation}
We assume that $\mathcal L$ is ortho-complemented. Let
$v\in{\LL^{[\perp \perp]}}$, with decomposition
\[
v=v_1+v_2,\quad v_1\in\LL,\quad{\rm and}\quad v_2\in {\LL^{[\perp]}}.
\]
Then, in view of \eqref{161112ineq}, $v_2=v-v_1\in{\LL^{[\perp
\perp]}}$, and so $v_2\in\LL^{[\perp]}\cap{\LL^{[\perp \perp]}}$.
Since
\[
\LL^{[\perp]}[+]{\LL^{[\perp \perp]}}=\V
\]
(recall that $\LL^{[\perp]}$ is also ortho-complemented; see
Proposition \ref{pn-ortho}), this implies that $v_2=0$ since $\V$
is non-degenerate. Thus there is equality in \eqref{161112ineq}.
\end{proof}

\section{Fundamental decompositions}
\setcounter{equation}{0}
A quaternionic inner product space ${\mathcal V}$ is
decomposable if it can be written as a direct and orthogonal sum
\begin{equation}
\label{decomposition1}
{\mathcal V}={\mathcal V}_+[\oplus]{\mathcal V}_-[\oplus]\mathcal N
\end{equation}
where ${\mathcal V}_+$ is a strictly positive subspace, ${\mathcal
V}_-$ is a strictly negative subspace, and $\mathcal N$ is a
neutral subspace. Representation \eqref{decomposition1} is called a fundamental
decomposition. A  quaternionic inner product space need not
be decomposable, and the decomposition will not be unique (unless
one of the spaces $\mathcal V_{\pm}$ is trivial).
A precise characterization of the decompositions
is given in the following results.

\begin{Pn}\label{4.4.1}
Assume that \eqref{decomposition1} holds. Then
$\mathcal N=\mathcal V\cap\mathcal V^{[\perp]}$ (that is $\mathcal N$ is equal to the isotropic part of $\mathcal
V$). \label{ch1-11.1-bognar}
\end{Pn}

\begin{proof} We first show that $\mathcal N\subset\mathcal V^0$.
Let $m\in\NN$, and let $v\in\V$ with decomposition
\begin{equation}
\label{decv}
v=v_++v_-+n,\quad\mbox{where}\quad v_\pm\in\V_\pm, \; n\in\NN.
\end{equation}
In view of \eqref{decomposition1} we have $[m,v_+]=[m,v_-]=0$. Furthermore, $[m,n]=0$
since the inner product vanishes in a neutral subspace (this is a
direct consequence of the Cauchy-Schwarz inequality). Thus $[m,v]=0$ and so
$m\in\V^0$.

\smallskip

Conversely, let $v_0\in\V^0$, with decomposition \eqref{decv}. Then,
\[
0=[v,v_+]=[v_+,v_+]
\]
and so $v_+=0$ since $\V_+$ is positive definite. Similarly,
$v_-=0$ and thus $v_0=n\in\mathcal N$.
\end{proof}
By  the definition of non-degenerate linear space, we have this immediate
consequence
of the previous result:
\begin{Cy}\label{Cor11.2}
All the decompositions of a  decomposable, non-degenerate inner product
space $\mathcal V$ are of the form
$$
\mathcal{V}=\mathcal{V}_+[\oplus]\mathcal{V}_-
$$
where $\mathcal{V}_+$ (resp. $\mathcal{V}_-$) is a strictly positive
(resp. negative) subspace.
\end{Cy}
The following is \cite[Lemma 11.4, p. 24-25]{bognar} in the
present setting.
\begin{Pn}\label{Prop62}
Let $\V$ be a  quaternionic non-degenerate inner product
space, and let $\LL$ be a positive definite subspace of $\LL$.
 There exists a fundamental decomposition of $\V$ with
$\mathcal V_+=\LL$ if and only if $\LL$ is maximal positive definite and
ortho-complemented.
\end{Pn}

\begin{proof} Assume first that
$\V=\LL[\oplus]\V_-[\oplus]\V^0$,
where $\V_-$ is negative definite and $\V^0$ is the isotropic
part of $\V$. Then $\LL$ is  ortho-complemented. Let $\mathcal M\supset\LL$
be a positive definite subspace containing $\LL$ and let
$v\in\mathcal M$, with decomposition
\[
v=v_++v_-+n,\quad v_+\in\LL,\quad v_-\in\LL_-,\quad n\in\V^0.
\]
By linearity, $v-v_+=v_-+n\in\mathcal M$. But
\[
[v-v_+,v-v_+]=[v_-,v_-]+[n,n]<0,
\]
unless $v_-=0$. But then $[v-v_+,v-v_+]=0$ implies $v=v_+$ (and
so $n=0$) since $\M$ is positive definite. Thus $v=v_+$ and
$\LL=\mathcal M$. Therefore,  $\LL$ is maximal positive definite.

\smallskip

Conversely, if $\LL$ is ortho-complemented, then
$\V=\LL[+]\LL^{[\perp]}$ and, since $\LL$ is positive definite,
the latter sum is direct, that is, $\V=\LL[\oplus]\LL^{[\perp]}$.
\smallskip

Since $\LL$ is maximal positive definite, it follows that
$\LL^{[\perp]}$ is negative.
Indeed, neither $\LL^{[\perp]}\setminus \LL$ nor  $\LL^{[\perp]}\cap
\LL$ contain positive
vectors $v$ since in the first case the space spanned by $v$ and
$\mathcal L$ would be positive,
contradicting the maximality of $\mathcal L$ and in the second case we
would have $[v,v]=0$
contradicting the positivity of $v$.
An application of Lemma
\ref{lemma:5.1} allows then to write $\LL^{[\perp]}$ as a direct
orthogonal sum of a negative definite space and of an isotropic
space $\NN$. Finally, the isotropic part $\NN$ of $\LL^{[\perp]}$
is the isotropic part of $\mathcal V$.
\end{proof}

%
To conclude this section we discuss some properties of linear
operators between  quaternionic inner product spaces.
The  linear operator $A$ will be called invertible if it is
one-to-one and its range is all of $\W$, or equivalently, if
there exists a  linear operator $B\,\,:\,\, \V\longrightarrow
\W$ such that $AB=I_{\W}$ and $BA=I_{\V}$.\\

Let ${\mathcal V}$ be a quaternionic inner product space which is
decomposable and non-degenerate, and let
\begin{equation}
\label{decomposition11}
{\mathcal V}={\mathcal
V}_+[\oplus]{\mathcal V}_-,
\end{equation}
where ${\mathcal V}_+$ is a strictly positive subspace and
${\mathcal V}_-$ is a strictly negative subspace. The map
\[
J(v)=v_+-v_-
\]
is called the associated fundamental symmetry. Since
$J(Jv)=v$, it follows that  $J$ is invertible, and
$J=J^{-1}$.  It is readily seen that
\begin{equation}
[v,w]=[Jv,Jw],\quad v,w\in \V.
\label{v=jv2}
\end{equation}
\begin{Tm}
\label{tm:preh} Let ${\mathcal V}$ be a decomposable and
non-degenerate quaternionic inner product space, and let
\eqref{decomposition1} be a fundamental decomposition of $\V$,
and let
\[
\langle v,w\rangle_J\stackrel{\rm def.}{=}[Jv,w],\quad v,w\in \V.
\]
Then,
\begin{equation}
\label{v=jv}
\langle v,w\rangle_J=[v,Jw]=[v_+,w_+]-[v_-,w_-],
\end{equation}
\begin{equation}
[v,w]=\langle v, Jw\rangle_J=\langle Jv,w\rangle_J,
\end{equation}
and $(\V,\langle\cdot, \cdot\rangle_J)$ is a pre-Hilbert space.
Furthermore, with $\|v\|_J=[ v,Jv]$
\begin{equation}
\label{eqcs}
|[v,w]|^2\le\|v\|_J^2\|w\|_J^2,\quad v,w\in \V.
\end{equation}
\end{Tm}
\begin{proof}
The first claim follows from the fact that both $\V_+$ and $\V_-$
are positive definite. In a quaternionic pre-Hilbert space, the
Cauchy-Schwarz inequality holds and this implies \eqref{eqcs}
since
\[
|[v,w]|^2=|\langle v,Jw\rangle_J|^2\le\|v\|_J^2\|Jw\|^2.
\]
Equations \eqref{v=jv} and \eqref{v=jv2} imply that
$\|w\|_J=\|Jw\|_J$, and this ends the proof.
\end{proof}
\begin{Rk}\label{Lemma11.2}
{\rm Let $\mathcal{V}$ be  a quaternionic, non-degenerate, inner product
vector space admitting a fundamental decomposition of the form
  $\mathcal{V}=\mathcal{V}_+[\oplus] \mathcal{V}_-$ and let $J$ be the
associated fundamental symmetry. Then $\mathcal{V}_+$ is $J$-orthogonal to
$\mathcal{V}_-$, i.e.
  $\langle v_+,w_-\rangle_J=0$ for every $v_+\in \mathcal{V}_+$ and
$w_-\in\mathcal{V}_-$, as one can see from formula (\ref{v=jv}).}
\end{Rk}

\section{Partial majorants}
\setcounter{equation}{0}

We now introduce and study some special topologies called partial
majorants. A standard reference for the material in this section
in the complex case is \cite[Chapter III]{bognar}. We begin by
proving a simple fact (which, in general, in not guaranteed in a
vector space over any field):
\begin{La}\label{lemma62}
Let $\mathcal{V}$ be a quaternionic inner product space and let
$w\in \mathcal{V}$. The maps
\begin{equation}
\label{sn3} v\mapsto\,\, p_w(v)=| [v,w]|, \quad v\in\mathcal{V}
\end{equation}
are semi-norms.
\end{La}
\begin{proof}
Property \eqref{sn1} is clear. Property \eqref{sn2} comes from
the fact  that the absolute value is multiplicative in $\mathbb
H$: $p_w(vc)=|[vc,w]|=|[v,w]c|=|[v,w]|\cdot |c|=|c| p_w(v)$.
\end{proof}

\begin{Dn}
The weak topology on $\mathcal{V}$ is the smallest topology such
that all the semi-norms \eqref{sn3} are continuous.
\end{Dn}

\begin{Dn}\label{partialmajorant}
$(a)$
A topology on the quaternionic indefinite inner product space
$\mathcal{V}$ is
called a partial majorant if it is locally convex and if all the
maps
\begin{equation}\label{maps}
v\quad \mapsto\quad[v,w]
\end{equation}
are continuous.\\
$(b)$ A partial majorant is called admissible if every continuous  linear function
from $\mathcal{V}$ to $\mathbb{H}$  is of the
form  $v\mapsto
[v,w_0]$ for some $w_0\in \mathcal{V}$.
\end{Dn}

\begin{Tm}
The weak topology of an inner
product space is a partial majorant. A locally convex topology is
a partial majorant if and only if it
is stronger than the weak topology.
\end{Tm}

\begin{proof}
To prove the first assertion, we have to show that in the weak topology
the maps (\ref{maps}) are continuous.
For any choice of $\varepsilon >0$, and for any $v_0,w\in\mathcal{V}$ the
inequality $|[v,w]- [v_0,w]|<\varepsilon$
is equivalent to $p_w(v-v_0)<\varepsilon$ and the set $\{v\in\mathcal{V}\ :\ p_w(v-v_0)<\varepsilon\}$ is a neighborhood
$U_{v_0}(p_w,\varepsilon)$ of $v_0$. Thus the weak topology is a partial majorant. \\
Let us now consider another locally convex topology stronger than the
weak topology. Then we have already shown that the inequality
$|[v,w]- [v_0,w]|<\varepsilon$ holds for $v\in
U_{v_0}(p_w,\varepsilon)$ which is also an open set in the
stronger topology and so any locally convex topology stronger than the weak topology is  a partial majorant. Finally, we consider a partial majorant. Let
$v_0, w_1,\ldots, w_n\in\mathcal{V}$, let $\varepsilon >0$. Then,
by definition, there are neighborhoods ${U}_\ell$ of $w_\ell$,
$\ell=1,\ldots, n$ such that for any $w\in {U}_{\ell}$ the
inequality $|[v,w_\ell]- [v_0,w_\ell]|<\varepsilon$, i.e.
$p_{w_\ell}(v-v_0)<\varepsilon$ holds. Thus any $w$ which belongs
to the neighborhood of $v_0$ given by $\cap_{\ell=1}^n {U}_\ell$
belongs to $U_{v_0}(p_{w_\ell}, \varepsilon)$ and the statement
follows.
\end{proof}
As a consequence we have:
\begin{Cy}
\label{corollary2.2p.85} Every partial majorant of a
non-degenerate inner product space $\V$ is Hausdorff.
\end{Cy}
\begin{proof}
Recall that any open set in the weak topology is
also open in the partial majorant topology. The weak topology is
Hausdorff if it separates points, i.e. if and only if for every
$w\in\mathcal{V}$ the condition $p_w(v)=|[v,w]|=0$ implies $v=0$.
But this is indeed the case since $\mathcal{V}$ is
non-degenerate.
\end{proof}

\begin{Pn}
\label{lemma2.4bognar}
%
If a topology is a partial majorant of the quaternionic inner product
space $\mathcal{V}$ then the orthogonal companion of every subspace is
closed.
\end{Pn}

\begin{proof}
Let $\mathcal{L}$ be a subspace of $\mathcal{V}$ and let
$\mathcal{L}^{[\perp]}$ its orthogonal companion. We show that
$\mathcal{L}^{[\perp]}$ is an open set. Let $v_0$ be in the
complement  $(\mathcal{L}^{[\perp]})^c$ of
$\mathcal{L}^{[\perp]}$; then there is $w\in \mathcal{L}$ such
that $[v_0,w]\not=0$. By continuity, there exists a neighborhood
${U}$ of $v_0$ such that $[v,v_0]\not=0$ for all $v\in{U}$, thus
$(\mathcal{L}^{[\perp]})^c$  is open.
\end{proof}
\begin{Cy}
\label{cyIII2.5}
If a topology is a partial majorant of a non-degenerate inner product
space
$\mathcal{V}$ then every ortho-complemented
subspace of $\mathcal{V}$ is closed.
\end{Cy}
\begin{proof}
Consider the subspace $\mathcal{L}^{[\perp]}$, orthogonal to
$\mathcal{L}$.
Then $\mathcal{L}^{[\perp\perp]}$ is closed by Proposition
\ref{lemma2.4bognar} and since $\mathcal{L}^{[\perp\perp]}=\mathcal{L}$ by
Proposition \ref{corollary9.5} the assertion follows.
\end{proof}
\begin{Cy}
Let $\tau$ be a partial majorant of the quaternionic inner product
$V$ and assume that $V$ is non-degenerate. Then the components of
any fundamental decompositions are closed with respect to $\tau$.
\label{cyIII.2.6}
\end{Cy}
\begin{proof}
This is a consequence of the previous corollary, since the two components
are orthocomplemented.
\end{proof}

\begin{Tm}\label{8.8.7} 
Let $\mathcal V$  be a non-degenerate quaternionic inner product
space and let $\tau_1$ and $\tau_2$ be two Fr\'echet partial
majorants of $\mathcal V$. Then, $\tau_1=\tau_2$.
\end{Tm}

\begin{proof}
Let $\tau$ be the topology $\tau_1\cup\tau_2$. Then we can show following
the proof of Theorem 3.3. p. 63 in \cite{bognar} that $\tau$ is a
Fr\'echet
topology stronger than $\tau_1$ and $\tau_2$.
We now consider the  two topological vector spaces $\mathcal{V}$ endowed
with $\tau$ and $\mathcal{V}$ endowed with $\tau_1$ and the identity map
acting between them.
By the closed graph theorem, see Theorem \ref{closedGR}, we have that the
identity map takes closed sets to closed sets and so $\tau_1$ is stronger
than $\tau$. A similar
argument holds by considering $\tau_2$ and thus $\tau=\tau_1=\tau_2$.
\end{proof}

Assume now that a partial majorant $\tau$ is defined by a norm $\|\cdot\|$
on a non-degenerate inner product space $\mathcal{V}$. Let us define
\begin{equation}
\label{polar}
\|v\|^\prime\stackrel{\rm def.}{=} \sup_{\|w\|\leq 1}|[v,w]|,
\qquad v\in\mathcal{V}.
\end{equation}
Then $\|\cdot \|'$ is a norm (called {\em polar} of the norm
$\|\cdot \|$), as it can be directly verified. As in the proof of
Lemma \ref{lemma62} the fact that the modulus is multiplicative in
$\mathbb{H}$ is what matters.
The topology $\tau'$ induced by $\|\cdot \|'$ is called the polar of the
topology $\tau$.\\
The definition (\ref{polar}) implies
\begin{equation}
\label{Buzenval-ligne-9}
|[v,w/\|w\|]|\leq\sup_{w\in\mathcal{V}} |[v,w/\|w\|]|\leq \sup_{ \|
w\|\leq 1}|[v,w]|=\| v\|',
\end{equation}
 from which we deduce the inequality  $|[v,w]|\leq \|v\|' \|w\|$. Thus the
polar of  a partial majorant is a partial majorant since (\ref{maps})
holds and thus one can define $\tau''\stackrel{\rm def.}{=}
(\tau')'$ and so on, iteratively.
\begin{Pn}
Let $\mathcal{V}$ be a non-degenerate inner product  space.
\begin{enumerate}
\item[(a)]
If $\tau_1$ and $\tau_2$ are normed partial majorants of $\mathcal{V}$. If
$\tau_1$ is weaker than $\tau_2$ then $\tau_2'$ is weaker than  $\tau_1'$.
\item[(b)] If $\tau$ be a normed partial majorant of $\mathcal{V}$, then
its polar $\tau'$ is a normed partial majorant on $\mathcal{V}$.
Furthermore, $\tau''\leq\tau$, and $\tau'''=\tau'$.
\end{enumerate}
\label{voltaire!!!}
\end{Pn}
\begin{proof} Let $\tau_1, \tau_2$ be induced by the norms $\|\cdot\|_1$
and $\|\cdot\|_2$, respectively
and let us assume that $\tau_1\leq \tau_2$. Then for
$w\in\mathcal{V}$ there exists $\lambda >0$ such that $\lambda \|
w\|_2\leq \| w\|_1$ and so, if we take $\|w\|_1\leq 1$ we have
$$
\sup_{\| w\|_1\leq 1}|[v,w]|\leq \sup_{\|w\|_2\leq 1}|[v,\lambda
w]|= \lambda \sup_{\|w\|_2\leq 1}|[v, w]|,
$$
so that $\tau_2'\leq\tau_1'$.\\
Moreover we have $\sup_{\|y\|'\leq 1}|[x,y]|\leq \| x\|$ and so
$\tau''\leq\tau$. Let us now use this inequality by replacing
$\tau$ by $\tau'$ and we get $\tau'''\leq \tau'$. By using point
$(a)$ applied to $\tau_1=\tau''$ and $\tau_2=\tau$ we obtain the
reverse inequality and so $\tau'''=\tau^\prime$.
\end{proof}
Among the partial majorants there are the admissible topology
(see Definition \ref{partialmajorant}). The next result shows that
an admissible topology which is metrizable is uniquely defined. In
order to prove the result, we recall that given a
quaternionic vector space $\mathcal{V}$, its conjugate
$\mathcal{V}^*$ is defined to be the  quaternionic vector
space in which the additive group coincides with $\mathcal{V}$
and whose multiplication by a scalar is given by $(c,v)\mapsto v\bar c
$. An inner product $(\cdot,\cdot)$ in $\mathcal{V}^*$ can be
assigned by $(v,w)\stackrel{\rm def.}{=}
[w,v]=\overline{[v,w]}$.
\begin{Tm} 
Let $\tau_1$, $\tau$ be admissible topologies on a quaternionic inner
product space $\mathcal{V}$. If $\tau_1$
is given by  a countable family of semi-norms, then $\tau_1$ is stronger
than $\tau$. Moreover, no more than one admissible topology of $\mathcal{V}$  is 
metrizable.
\label{Mairie-de-Montreuil-ligne-9}
\end{Tm}
\begin{proof}
Assume that $\tau_1$ and $\tau$ are given by the families of semi-norms
$\{p_i\}$, $i\in\mathbb{N}$, and $\{q_\gamma\}$, $\gamma\in\Gamma$,
respectively.
By absurd, suppose that $\tau_1$ is not stronger than $\tau$. Then there
exists an open set in $\tau$ that does not contain any open set in
$\tau_1$ and, in particular, it does not contain
$$
\{ v\in\mathcal{V} \ |\  p_i(v)<\frac 1n, \ \ i=1,\ldots ,n, \ {\rm for}\,
n\in \mathbb{N} \}.
$$
Thus, there exists a sequence $\{v_n\}\subset\mathcal{V}$ such that
$p_i(v_n)<\frac 1n$ but $\max_{k=1,\ldots
,m}q_{\gamma_k}(v_n)=q_{\gamma_j}(v_n)\geq\varepsilon$ for some
$\varepsilon>0$. By choosing $w_n=nv_n$ we have
\begin{equation}\label{58}
\max_{i=1,\ldots, n} p_i(w_n)< 1, \qquad q_{\gamma_j}(w_n)\geq
n\varepsilon, \quad n\in\mathbb{N}.
\end{equation}
Let us consider the subspace of $\mathcal{V}$ given by
$\mathcal{L}=\{v\in\mathcal{V}\ |\ q_{\gamma_j}(v)=0\}$ and the quotient
$\hat{\mathcal{L}}\stackrel{\rm def.}{=}
\mathcal{V}/\mathcal{L}$. We can endow $\hat{\mathcal{L}}$ with the norm
$\| \hat v\|\stackrel{\rm def.}{=}
q_{\gamma_j}(v)$, for $\hat {v}=v+\mathcal{L}\in\hat{\mathcal{L}}$. Let
$\hat{\varphi}:\ \hat{\mathcal{L}}\to\mathbb{H}$ be a linear  function
which is also continuous (bounded):
$$
|\hat{\varphi}(\hat{v})|\leq \|\hat{\varphi}\|\,
\|\hat{v}\|,\quad\hat{v}\in\hat{\mathcal{L}}.
$$
Then the formula $\varphi(v)\stackrel{\rm def.}{=}\hat{\varphi}(\hat{v})$,
$v\in\mathcal{V}$, $v\in\hat{v}$, defines a linear and
continuous function on $\mathcal{V}$ since
$$
|\varphi(v)|\leq \|\hat{\varphi}\| \, \|\hat{v}\|= \|\hat{\varphi}\|
q_{\gamma_j}(v).
$$
Thus $\varphi$ is continuous in the topology $\tau$ and since $\tau$ is
admissible, $\varphi(v)=[v,w_0]$ for some suitable $w_0\in\mathcal{V}$.
We conclude that $\varphi$ is also continuous in the topology $\tau_1$. So
for some  $r\in\mathbb{N}$ and $\delta >0$ we have
$$
|\varphi(v)| \leq \frac{1}{\delta}\max_{i=1,\ldots ,r} p_i(v), \quad
v\in\mathcal{V}.
$$
This last inequality together with (\ref{58}) give
$|\varphi(w_n)|<1/\delta$ for $n>r$. So the sequence
$\{\hat{\varphi}(\hat{w}_n)\}$ is bounded
for any $\hat{\varphi}$ fixed in the conjugate space $\hat{\mathcal{L}}^*$
of the normed space $\hat{\mathcal{L}}$. However, we can look at
$\hat{\varphi}(\hat{w}_n)$
as the value of the functional $\hat{w}_n$ acting on the elements of the
Banach space $\hat{\mathcal{L}}^*$.
 Since we required that
$|\hat{\varphi}(\hat{v})|\leq \|\hat{\varphi}\|\, \|\hat{v}\|$,
for $\hat{v}\in\hat{\mathcal{L}}$ the functional $\hat{w}_n$ is
continuous. By the quaternionic version of the Hahn-Banach
theorem, see e.g. \cite[Theorem 4.10.1]{MR2752913}, we deduce
that $\|\hat{w}_n\|=q_{\gamma_j}(\hat{w}_n)$. From (\ref{58}),
more precisely from $q_{\gamma_j}(w_n)\geq n\varepsilon$, we
obtain a contradiction with the principle of uniform boundedness,
see Theorem \ref{tmunivb}.
\end{proof}

\section{Majorant topologies and inner product spaces}
\setcounter{equation}{0}

The material in this section can be
found, in the complex case,  in  \cite[Chapter IV]{bognar}.
\begin{Dn}
A locally convex topology on $(\mathcal{V},[\cdot,\cdot])$ is
called a majorant if the inner product is jointly continuous in
this topology. It is called a complete majorant if it is
metrizable and complete. It is called a normed majorant if it is
defined by a single (semi-)norm, and a Banach majorant if it is
moreover complete with respect to this norm. It is called a
Hilbert majorant if it is a complete normed majorant, and the
underlying norm is defined by an inner product. \label{def:maj}
\end{Dn}

Of course, the norm defining a Banach majorant (and hence the
inner product defining a Hilbert majorant) is not unique. But it
follows from Theorem \ref{OPNMAP} that any two such norms are
equivalent.

\begin{Pn}\mbox{}\\
$(a)$ Given a majorant, there exists a weaker majorant defined by
a single semi-norm.\\
$(b)$ A normed partial majorant $\tau$ on the non-degenerate inner
product space $\mathcal V$ is a majorant if and only if it is
stronger than its polar: $\tau^\prime \le \tau$.
\label{nation-ligne-9}
\end{Pn}

\begin{proof}
$(a)$ From the definition of a majorant, there exist semi-norms
$p_1,\ldots , p_N$ and $\epsilon>0$ such that
\[
|[u,v]|\le 1, \quad\forall u, v\in {U},
\]
where
\[
{U}=\left\{v\in\mathcal V\,;\,p_j(v)\le\epsilon,\,\,j=1,\ldots N
\right\}.
\]
It follows that the inner product is jointly continuous with
respect to the semi-norm $\max_{j=1,\ldots N}p_j$.\\

$(b)$ Recall that the polar $\tau^\prime$ is defined by
\eqref{polar}. We have $\tau^\prime \le \tau$ if and only if the
identity map from $(\mathcal V,\tau)$ into $(\mathcal
V,\tau^\prime)$ is continuous, that is if and only if there
exists $k>0$ such that
\begin{equation}
\|v\|^\prime\le k\|v\|,\quad \forall v\in\mathcal V.
\end{equation}
This is turn holds if and only if
\begin{equation}
|[v,u]|\le k\|v\|,\quad \forall v,u\in\mathcal V\quad{\rm with}\quad
\|u\|\le 1.
\end{equation}
The result follows since  any such $w\not=0$ is such that $\|w\|\le 1$ if
and only if it be written as $\frac{w}{\|w\|}$,
for some $w\not=0\in\mathcal V$.
\end{proof}

\begin{Pn} 
Let $\mathcal V$ be a non-degenerate inner product space, admitting
a normed majorant. Then there exists a
weaker normed majorant which is self-polar.
\label{avron-ligne-9}
\end{Pn}

\begin{proof} We briefly recall the proof of \cite[p. 85]{bognar}. The key
is that the polar norm (defined in
\eqref{polar})
is still a norm in the quaternionic case. By maybe renormalizing we assume
that
\begin{equation}
\label{Maraichers-ligne-9}
|[u,v]|\le \|u\|\|v\|,\quad u,v\in\V,
\end{equation}
where $\|\cdot\|$ denotes a norm defining majorant. Define a sequence of
norms $(\|\cdot\|_n)_{n\in\mathbb N}$ by
$\|\cdot\|_1=\|\cdot\|$ and
\begin{equation}
\label{ender-s-game}
\|u\|_{n+1}=\left(\frac{1}{2}(\|u\|_n^2+(\|u\|_n^\prime)^2)\right)^{\frac{1}{2}},\quad
n=1,2,\ldots,
\end{equation}
where we recall that $\|\cdot\|^\prime$ denotes the polar norm of
$\|\cdot\|$; see \eqref{polar}. An induction shows that
each $\|\cdot\|_n$ satisfies \eqref{Maraichers-ligne-9} and that the
sequence $(\|\cdot\|_n)_{n\in\mathbb N}$  is decreasing,
and thus defining a semi-norm
$\|\cdot\|_\infty=\lim_{n\rightarrow\infty}\|\cdot\|_n$. One readily shows
that
$\|\cdot\|_\infty\ge \frac{1}{\sqrt{2}}\|\cdot\|_1^\prime$, and hence
$\|\cdot\|_\infty$ is a norm, and a majorant since
it also satisfies \eqref{Maraichers-ligne-9} by passing to the limit the
corresponding inequality for $\|\cdot\|_n$.

\smallskip

We now show that the topology defined by $\|\cdot\|_\infty$ is self-polar.
We first note that the sequence of polars
$(\|\cdot\|^\prime_n)_{n\in\mathbb N}$  is increasing, and bounded by the
polar $\|\cdot\|_\infty^\prime$. Set
$\|\cdot\|_e=\lim_{n\rightarrow\infty}\|\cdot\|_n^\prime$.
Applying inequality \eqref{Buzenval-ligne-9} to $\|\cdot\|_n$ and taking
limits leads to
\[
|[u,v]|\le \|u\|_e\|v\|_\infty,\quad u,v\in\V.
\]
Thus $\|\cdot\|_\infty^\prime\le\|\cdot\|_e$, and we get that
$\|\cdot\|_\infty=\|\cdot\|_e$. Letting $n\rightarrow\infty$ in
\eqref{ender-s-game} we get
$\|\cdot\|_\infty=\|\cdot\|_\infty^\prime$.
\end{proof}

\begin{Pn}
\label{tm4.2p86}
Let $(\mathcal V,[\cdot,\cdot])$ be a quaternionic non-degenerate
inner product space. Then a partial majorant is a minimal
majorant if and only if it is normed and self-polar.
\end{Pn}

\begin{proof} Assume first that the given partial majorant $\tau$ is a
minimal majorant.
By item $(a)$ of Proposition \ref{nation-ligne-9} there is a
weaker majorant $\tau_a$  defined by a single semi-norm. Moreover
by Corollary \ref{corollary2.2p.85} any partial majorant
(and in particular any majorant) is Hausdorff, and so the
$\tau_a$ is Hausdorff and the above semi-norm is in fact a norm. By
Proposition \ref{avron-ligne-9} there exists a
self-polar majorant $\tau_\infty$ which is weaker that $\tau_1$.  The
minimality of $\tau$ implies that $\tau_\infty=\tau$.

\smallskip

Conversely, assume that the given partial majorant $\tau$ is
normed and self-polar. Then $\tau$ is a majorant in view of item
$(b)$ of Proposition \ref{nation-ligne-9}. Assume that
$\tau_a\le\tau$ is another majorant. Then, by part $(b)$ in Lemma
\ref{nation-ligne-9}, $\tau_a\ge \tau_a^\prime$,
and by item $(a)$ of Proposition \ref{voltaire!!!} we have
$\tau_a^\prime\ge \tau^\prime$. This ends the proof since
$\tau$ is self-polar.
\end{proof}

\begin{Tm}\label{9.9.5}
Let $\mathcal V$ be a quaternionic non-degenerate
inner product space, and let $\tau$ be an admissible topology
which is moreover a majorant. Then $\tau$ is minimal, it defines a
Banach topology and is the unique admissible majorant on
$\mathcal V$. Finally, $\tau$ is stronger than any other admissible
topology on $\mathcal V$.
\end{Tm}

We now introduce the Gram operator. It will play an important role in the
sequel. Recall that Hilbert majorants have been
defined in Definition \ref{def:maj}.

\begin{Pn}\label{9.9.6}
Let $(\mathcal V,[\cdot,\cdot])$ be a quaternionic inner-product
space, admitting a Hilbert majorant, with associated inner
product $\langle\cdot, \cdot\rangle$, and corresponding norm
$\|\cdot\|$. There exists a linear continuous operator $G$,
self-adjoint with respect to the inner product $\langle\cdot,
\cdot\rangle$, and such that
\[
\label{def:gram}
[v,w]=\langle v, Gw\rangle,\quad v,w\in\mathcal V.
\]
\end{Pn}

\begin{proof} The existence of $G$ follows from Riesz' representation
for continuous functionals, which still holds in quaternionic
Hilbert spaces (see \cite[p. 36]{bds}, \cite[Theorem II.1, p.
440]{MR768240}); the fact that $G$ is Hermitian follows from the
fact that the form $[\cdot,\cdot]$ is Hermitian. In the complex
case, an everywhere defined Hermitian operator in a Hilbert space
is automatically bounded; rather than proving the counterpart of
this fact in the quaternionic setting we note, as in  \cite[p.
88]{bognar} that there exists a constant $k$ such that
\begin{equation}
\label{eq:maj}
|[u,w]|\le k\|u\|\cdot\|v\|,\quad \forall u,v\in\mathcal V.
\end{equation}
The boundedness of $G$ follows from \eqref{eq:maj} and
$[v,Gv]=\|Gv\|^2$.
\end{proof}

The semi-norm
\begin{equation}
v\mapsto\|Gv\|
\label{Croix-de-Chavaux-ligne-9}
\end{equation}

 defines a topology called the Mackey topology. As we remarked after
Definition
\ref{def:maj} the inner product defining a given Hilbert majorant is not
unique, and so to every inner product will
correspond a different Gram operator.
\begin{Pn}\label{9.9.7}
The Mackey topology is admissible and is independent of the choice of the
inner product defining the Hilbert majorant.
\end{Pn}
\begin{proof} The uniqueness will follow from Theorem
\ref{Mairie-de-Montreuil-ligne-9} once we know that the topology,
say $\tau_G$,
associated to the semi-norm \eqref{Croix-de-Chavaux-ligne-9} is
admissible. From the inequality
\[
|[u,v]|=\langle Gu,v\rangle\le\|Gu\|\cdot\|v\|
\]
we see that $\tau_G$ is a partial majorant. To show that it is admissible,
consider a linear functional $f$ continuous with
respect to $\tau_G$. There exists $k>0$ such that
\[
|f(u)|\le k\|Gu\|,\quad \forall u\in\V.
\]
Thus the linear relation
\[
(kGu,f(u)),\quad u\in\V
\]
is the graph of a contraction, say $T$,
\[
T(Gu)=\frac{1}{k}f(u),\quad\forall u\in\V,
\]
in the pre-Hilbert space $({\rm Ran}~G)\times \H$, the latter being
endowed with the inner product
\[
\langle (Gu,p), (Gv,q)\rangle_{\V\times \mathbb H}=\langle
Gu,Gv\rangle+\overline{q}p=[Gu,v]+\overline{q}p.
\]
The operator $T$ admits a contractive extension to all of
$\V\times \H$, and by Riesz representation theorem, there exists
$f_0\in\V$ such that
\[
T(u)=\langle u, f_0\rangle,\quad \forall u\in \mathcal V.
\]
Thus
\[
f(u)=kT(Gu)=k\langle Gu,f_0\rangle=[u,kf_0],
\]
which ends the proof.
\end{proof}

Consider a subspace $\mathcal L$ of a  quaternionic
inner-product space $(\mathcal V,[\cdot,\cdot])$, the latter
admitting an Hilbert majorant with associated inner product
$\langle\cdot,\cdot\rangle$ and associated norm $\|\cdot\|$. We
denote by $P_{\mathcal L}$ the orthogonal projection onto
$\mathcal L$ in the Hilbert space
$(\V,\langle\cdot,\cdot\rangle)$, and set
\begin{equation}
G_{\mathcal L}=P_{\mathcal L}G\big|_{\mathcal L}.
\label{Robespierre-ligne-9}
\end{equation}
\begin{Pn}\label{9.9.8}
Consider $\mathcal V$ be a quaternionic
inner-product space, admitting an Hilbert majorant, let $\mathcal
L$ be a closed subspace of $\mathcal V$ and let $G_{\mathcal L}$
be defined by \eqref{Robespierre-ligne-9}. Then:\\
$(a)$ An element $v\in\V$ admits a projection onto $\mathcal L$ if and
only if
\begin{equation}
\label{le-trocadero}
P_{\mathcal L}v\in{\rm ran}~G_{\mathcal L}.
\end{equation}
$(b)$ $\mathcal L$ is ortho-complemented in $(\mathcal
V,[\cdot,\cdot])$ if and only if
\[
{\rm ran}~P_{\mathcal L}G={\rm ran}~G_{\mathcal L}.
\]
\end{Pn}

\begin{proof}
$(a)$ The vector $v\in V$ has a (not necessarily unique)
projection, say $w$ on $\LL$ if and only if
\[
[v-w,u]=0,\quad \forall u\in \LL,
\]
that is, if and only if
\[
\langle G(v-w),u\rangle=0,\quad \forall u\in \LL.
\]
This last condition is equivalent to
$P_{\LL} Gv=G_{\LL} w$,
which is equivalent to \eqref{le-trocadero}.

\smallskip

$(b)$ The second claim is equivalent to the fact that every
element admits a projection on $\LL$, and therefore follows from
$(a)$.
\end{proof}

\section{The spectral theorem and decomposability}
\setcounter{equation}{0}
The spectral theorem for Hermitian operators is stated in
\cite{MR0137500},
\cite{viswanath_thesis},
\cite{MR44:2067}  in which, however,
a proof is not provided. Moreover, in these works, the spectrum used is
not the
$S$-spectrum, see \cite[p. 141]{MR2752913}, thus for the sake of completeness
we state and prove the result. To this end, we need some preliminaries.

\smallskip

We first note that any  linear quaternionic
Hilbert space $\mathcal V$ can be also considered as a complex
Hilbert space, its so-called symplectic image denoted by
$\mathcal V_s$, which coincides with $\mathcal V$ as Abelian
additive group and  whose multiplication by a scalar is the
multiplication given in $\mathcal V$  restricted to $\mathbb C$.
Here we identify $\mathbb C$ with the set of
quaternions of the form $x_0+ix_1$. Any linear operator $T$ on
$\mathcal V$ is obviously also $\mathbb C$-linear and so it is a
linear operator on $\mathcal V_s$. We denote by $T_s$ the
operator $T$ when it acts on $\mathcal V_s$.  The converse is not
true, i.e. if $S$ is a $\mathbb C$-linear operator
acting on $\mathcal V_s$ then $S$ is not, in general, a
linear operator on $\mathcal V$,
unless additional hypothesis are given.
It is immediate to verify that if $T$ is Hermitian then $T_s$ is
Hermitian (see also \cite{viswanath_thesis}).\\
We now state the spectral theorem:
\begin{Tm}
Let $A$ be a Hermitian operator on the  quaternionic Hilbert
space $\mathcal V$. Then there exists a spectral measure $E$
defined on the Borel sets in $\mathbb R$ such that
\begin{equation}
\label{spectraltheo}
A=\int_{-\infty}^{+\infty} \lambda dE (\lambda) .
\end{equation}
\end{Tm}
\begin{proof}
We observe that if $A$ is a Hermitian  linear operator, then its
$S$-spectrum is real.\\
Then we consider the symplectic image  $\mathcal V_s$ of
$\mathcal V$ and the operator $A_s$ which is Hermitian and whose
(real) spectrum coincide with the spectrum of $A$. Then we can
use the classical spectral theorem to write
$$
A_s=\int_{-\infty}^{+\infty} \lambda dE_s (\lambda)
$$
where $dE_s (\lambda)$ is a spectral measure with values in the
lattice of projections in $\mathcal V_s$. Since the support of
$E$ is contained in $\mathbb R$, we use Corollary 6.1 in
\cite{viswanath_thesis} to guarantee that $E$ is a spectral
measure with values in the lattice of projections in $\mathcal
V$. This concludes the proof since $A_s$ is in fact $A$.
\end{proof}

\begin{Tm}\label{10.10.2}
Let $(\mathcal V,[\cdot,\cdot])$ be a quaternionic inner-product
space, admitting a Hilbert majorant. Then $\mathcal V$ is
decomposable, and there exists a fundamental decomposition such
that all three components and any sum of two of them are complete
with respect to the Hilbert majorant.
\end{Tm}

\begin{proof} 
As in the proof of the corresponding result in the complex case (see
\cite[p. 89]{bognar} we apply
the spectral theorem to the Gram operator $G$ associated to the form
$[\cdot, \cdot]$, and write $G$ as
\eqref{spectraltheo}:\
\[
G=\int_{-\infty}^{+\infty} \lambda dE (\lambda),
\]
where the spectral measure is  continuous and its support is
finite since $G$ is bounded. We then set
\[
\V_-=E(0^-)\V,\quad \V_0=(E(0)-E(0^-))\V,\quad{\rm and}\quad
\V_+=(I-E(0))\V.
\]
We have
\[
\V=\V_-[\oplus]\V_0[\oplus]\V_+.
\]
Each of the components and each sum of pairs of components of
this decomposition is an orthogonal companion, and therefore
closed for the Hilbert majorant in view of Proposition
\ref{lemma2.4bognar}.
\end{proof}

In the next result, the space is non-degenerate, but the
majorant is a Banach majorant rather than a Hilbert majorant.

\begin{Pn}\label{10.10.3}
Let $(\mathcal V,[\cdot,\cdot])$ be a quaternionic non-degenerate
inner-product space, admitting a Banach majorant $\tau$ and a
decomposition majorant $\tau_1$. Then, $\tau_1\le \tau$.
\end{Pn}
\begin{proof}
Let $\mathcal V=\mathcal V_+[\oplus]\mathcal V_-$ be a fundamental
decomposition of $\mathcal{V}$. By Corollary \ref{cyIII.2.6} the space
$\mathcal V_+$ is closed in the topology $\tau$. Let $P_+$ denote the map
\begin{equation}\label{P+}
P_+v=v_+
\end{equation}
where $v=v_++v_-$ is the decomposition of $v\in\mathcal V$ along
the given fundamental decomposition of $\mathcal V$. We claim
that the graph of $P_+$ is closed, when $\mathcal V$ is endowed
with the topology $\tau$. Indeed, if $(v_n)_{n\in\mathbb N}$ is a
sequence converging (in the topology $\tau$) to $v\in\mathcal V$ and
such that the sequence $((v_n)_+)_{n\in\mathbb N}$ converges to
$z\in\mathcal V_+$ also in the topology $\tau$. Since the inner
product is continuous with respect to $\tau$ we have for
$w\in\mathcal V_+$
\begin{align*}
[z-v_+,w]&=\lim_{n\rightarrow\infty}[(v_n)_+,w]-[v_+,w]\\
&=\lim_{n\rightarrow\infty}[v_n,w]-[v_+,w]
=[v,w]-[v_+,w] =[v-v_+,w]=0
\end{align*}
and so $z=v_+$. By the closed graph theorem (see Theorem
\ref{closedGR}) $P_+$ is continuous. The same holds for the
operator $P_-v=v_-$ and so the operator
\begin{equation}
\label{eqJ}
Jv=v_+-v_-
\end{equation}
is continuous from $(\mathcal V,\tau)$ onto $(\mathcal V,\tau)$.
Recall now that $[Jv,v]$ is the square of the $J$-norm defining
$\tau_1$. We have
\[
[Jv,v]\le k\|Jv\|\cdot\|v|,
\]
where $\|\cdot\|$ denotes a norm defining
$\tau$,  and
\[
[Jv,v]\le k\|Jv\|\cdot\|v|\le k_1\|v\|^2
\]
since $J$ is continuous.
It follows that the inclusion map is continuous from $(\mathcal
V,\tau)$ into $(\mathcal V,\tau_1)$, and so $\tau_1\le \tau$.
\end{proof}

\begin{Pn}\label{10.10.4} 
Every decomposition majorant is a minimal majorant.
\label{10104}
\end{Pn}

\begin{proof}
A decomposition majorant is in particular a partial majorant and is normed
(with associated $J$-norm
$\|u\|_J=[Ju,u]$, where $J$ is associated to the decomposition as in
\eqref{eqJ}). Thus, using
Proposition \ref{tm4.2p86}, to prove the minimality it is enough to
show that $\|u\|_J$ is self-polar. That this holds follows from
\[
\|u\|_J^\prime=\sup_{\|v\|_J\le 1}|[Ju,v]|=\sup_{\|v\|_J\le
1}|[u,Jv]|=\|u\|_J.
\]
\end{proof}

The question of uniqueness of a minimal majorant is considered in the next
proposition.

\begin{Pn} \label{10.10.6}
Let $(\mathcal V,[\cdot,\cdot])$ be a quaternionic inner-product
space, admitting a decomposition
\begin{equation}
\label{porte-de-montreuil-ligne-9}
\V=\V_+[\oplus]\V_-,
\end{equation}
where $\V_+$ is positive definite and $\V_-$ is negative definite. Assume
that $\mathcal{V}_+$ (resp. $\mathcal{V}_-$)
is intrinsically complete. Then,
so is $\V_-$ (resp. $\V_+$). Then $(\mathcal V,[\cdot,\cdot])$ has a
unique minimal majorant.
\label{pn6.2}
\end{Pn}

\begin{proof}
The topology $\tau$ defines a fundamental decomposition, and an associated
minimal majorant $\|\cdot\|_J$.
See Proposition \ref{10104}.
Let $\tau$ be another minimal majorant. By Proposition \ref{tm4.2p86} it
is normed and self-polar and so there
is a norm $\|\cdot\|$ and $k_1>0$ such that
\[
\|v_+\|\le k_1\sup_{\substack{y\in\V_+\\ \|v\|\le 1}}|[v_+,y]|.
\]
Using the uniform boundedness we find $k_2>0$ such that
\[
|[v_+,y]|\le k_2[v_+,v_+],\quad \forall y\,\,\mbox{\rm such that}\,\,
\|y\|\le 1.
\]
Hence, with $C=k_1k_2$,
\begin{equation}
\label{mnb}
\|v_+\|\le C[v_+,v_+],\quad \forall v_+\in\V_+.
\end{equation}
Let now $v\in\V$ with decomposition $v=v_++v_-$, where $v_\pm\in\V_\pm$.
Since $\tau$ is a normed
majorant, there exists $C_1$ such that
\[
\|v_+\|2\le C[v_+,v_+]=C[v_+,v]\le CC_1\|v_+\|\cdot\|v\|
\]
Hence
\[
\|v\|_J^2=[Jv,v]\le C_1\|Jv\|\cdot\|v\|
=C_1\|2v_+-v\|\cdot K\|v\|^2
\]
for an appropriate $K>0$. The identity map is there continuous from
$(\mathcal{V},\tau)$ onto $(\mathcal{V},\|\cdot\|_J)$. Since $\tau$ is
defined by a single norm, it follows that the identity map is also
continuous from $(\mathcal{V},\|\cdot\|_J)$
onto $(\mathcal{V},\tau)$ and this ends the proof.
\end{proof}

\begin{Pn}\label{10.10.7}
Let $(\mathcal V,[\cdot,\cdot])$ be a quaternionic inner-product
space, admitting a decomposition of the form
\eqref{porte-de-montreuil-ligne-9}, and with associated
fundamental symmetry $J$. Then: $(a)$ Let $\mathcal L$ denote
a positive subspace of $\V$. Then, the operator $P_+\big|_{\mathcal L}$
and its inverse are $\tau_J$ continuous.\\
$(b)$ Given another decomposition of the form
\eqref{porte-de-montreuil-ligne-9}, the positive (resp. negative)
components are simultaneously intrinsically complete.
\end{Pn}

\begin{proof}
To prove the result we follow \cite[pp. 93-94]{bognar}.
Let $\mathcal L$ be a positive subspace of $\mathcal V$ and let
$v\in\mathcal L$.
By recalling (\ref{v=jv}), (\ref{P+}) and setting $P_-v=v_-$, where
$v=v_++v_-$ is the decomposition of $v$ with respect to the fundamental
decomposition $\mathcal V=\mathcal V_+[\oplus] \mathcal V_-$, we have:
\[
\| v \|_J^2=\|P_+v\|_J^2+\|P_-v\|_J^2.
\]
Since $\mathcal V_+$ and $\mathcal V_-$ are $J$-orthogonal, see Remark
\ref{Lemma11.2}, we then have
$$
[v,v]=\|P_+v\|_J^2-\|P_-v\|_J^2
$$
and so, since $\mathcal L$ is positive,
\[
\| v \|_J^2=2\|P_+v\|_J^2- [v,v] \le 2\|P_+v\|_J^2.
\]
It is immediate that $\| P_+ v \|_J^2\le\| v \|_J^2$ and so we conclude
that both $P_+$ and its inverse are $\tau_J$ continuous as stated in
(a).\\
To show $(b)$, we assume that there is another fundamental decomposition
$\mathcal V=\mathcal V'_+[\oplus] \mathcal V'_-$. If we suppose that
$\mathcal{V}'_+$ is intrinsically complete, then Proposition \ref{10.10.6}
implies that $\mathcal V'_+$ is complete with respect to the decomposition
majorant
corresponding to the decomposition $\mathcal V=\mathcal V_+[\oplus]
\mathcal V_-$. Part (a) of the statement implies that also $P^+\mathcal
V'_+$ is complete in this topology and so it is intrinsically complete. If
$P^+ \mathcal V'_+=\mathcal V_+$ there is nothing to prove.
Otherwise there exists a non-zero $\tilde v\in\mathcal V_+$
orthogonal to $P^+ \mathcal V'_+$ so $\tilde v$ is orthogonal to
$\mathcal V'_+$. Then the subspace
$\mathcal U$ spanned by $\tilde v$
and $\mathcal V'_+$ is positive. Indeed, for a generic nonzero element
$u=\tilde v+\tilde v^\prime$
($\tilde v^\prime\in\mathcal V'_+$)
we have
$$
[u,u]=[\tilde v+\tilde v^\prime, \, \tilde v+\tilde v^\prime]=[\tilde v,
\, \tilde v]+[\tilde v^\prime, \,
\tilde v^\prime]>0.
$$
This
implies that $\mathcal U$ is a proper extension of $\mathcal V'_+$ which
is absurd by Proposition \ref{Prop62}. This completes the
proof.\end{proof}

\section{Quaternionic Krein spaces}
\setcounter{equation}{0}
In this section we will study quaternionic Krein spaces following
\cite[Chapter V]{bognar}. As in the classical case, they are characterized
by the fact that they are inner product spaces non-degenerate,
decomposable and complete. We will show that the scalar product associated
to the decomposition gives a norm, and so a topology, which does not
depend on the chosen decomposition. We will also study ortho-complemented
subspaces of a Krein space and we will prove that they are closed
subspaces which are Krein spaces themselves.
\begin{Dn}
If a quaternionic inner product space ${\mathcal K}$ has a fundamental
decomposition
\begin{equation}\label{krein}
{\mathcal K}={\mathcal K}_+[\oplus]{\mathcal K}_-,
\end{equation}
where ${\mathcal K}_+$ is a strictly positive subspace while ${\mathcal
K}_-$
is strictly negative and if ${\mathcal K}_+$ and ${\mathcal K}_-$ are intrinsically complete, then
we say that $\mathcal K$ is a {\em Krein space}.
\end{Dn}
The decomposition of a Krein space is obviously not unique when one of the
components is not trivial.
Both the spaces $\mathcal{K}_+$ and $\mathcal{K}_-$ are Hilbert spaces and
they can be, in particular,
of finite dimension. The Krein space is then called a Pontryagin space
when $\V_-$ is finite dimensional.

\begin{Pn}\label{11.11.2} 
A Krein space is non-degenerate and decomposable. Each fundamental decomposition has  intrinsically complete components $\mathcal K_{\pm}$.
\end{Pn}
\begin{proof}
A Krein space is obviously decomposable by its definition and
non-degenerate by Proposition \ref{4.4.1}.\\
By Theorem \ref{10.10.7} $(b)$, given (\ref{krein}) and any other
fundamental decomposition
${\mathcal K}={\mathcal K}'_+[\oplus]{\mathcal K}'_-$ if ${\mathcal K}_+$
is intrinsically complete so is ${\mathcal K}'_+$ (and similarly for
${\mathcal K}'_-$).
\end{proof}
\begin{Pn} \label{11.11.3} 
A non-degenerate, decomposable, quaternionic inner product space $\mathcal
K$ is a Krein space if and only if for every associated fundamental
symmetry
$J$, $\mathcal K$ endowed with the inner product $\langle v,
w\rangle_J\stackrel{\rm def.}{=}[v,Jw]$ is  a Hilbert space.
\end{Pn}
\begin{proof}
Let $\mathcal K$ be a non-degenerate, decomposable, quaternionic inner
product space, i.e. $\mathcal{K}=\mathcal{K}'_+[\oplus] \mathcal{K}'_-$.
If $\mathcal K$ is a Krein space then the associated fundamental symmetry
$J=P^+-P^-$ makes it into a pre-Hilbert space, see Theorem \ref{tm:preh}.
Completeness follows from the fact that $\mathcal{K}_\pm$ are both
complete.
Conversely, assume that given a fundamental symmetry
$J$ the inner product $\langle v, w\rangle_J\stackrel{\rm def.}{=}
[v,Jw]$ makes $\mathcal K$ a Hilbert space. The intrinsic norm in
$\mathcal{K}_+$ is obtained by restricting the $J$-inner product to
$\mathcal{K}_+$. Any Cauchy sequence in $\mathcal K_+$ converges to an
element in $\mathcal K$ and it is immediate to verify that this element
belongs to $\mathcal K_+$.
\end{proof}
\begin{Tm}\label{Tm:V1.3}
Let $\mathcal K$ be a  quaternionic vector space with inner product
$[\cdot,\cdot]$. Then $\mathcal K$ is a Krein space if and only if :
\begin{enumerate}
\item[(a)] $[\cdot,\cdot]$ has a Hilbert majorant $\tau$ with associated
inner product $\langle\cdot,\cdot\rangle$ and norm $\|v\|=\langle
v,v\rangle^{\frac 12}$;
\item[(b)]  the Gram operator of $[\cdot, \cdot]$ w.r.t. $\langle\cdot,\cdot\rangle$, i.e., the operator
$G$ which satisfies $[v,w] = \langle v,Gw\rangle$, $v,w \in\mathcal{K}$, is boundedly invertible.
\end{enumerate}
\end{Tm}
\begin{proof} We follow the proof of Theorem V, 1.3 in \cite{bognar}, by
repeating the main arguments.
Assume that $\mathcal{K}$ is a Krein space and denote by $J$ the
fundamental symmetry associated to the chosen decomposition (\ref{krein}).
Define a norm using the $J$-inner product $\langle\cdot ,\cdot\rangle_J$
and let $\tau_J$ be the corresponding topology which is is a decomposition
majorant
by Proposition \ref{10.10.6} and a Hilbert majorant. Since
\[
[v,w]= [v, J^2w]=\langle v, Jw\rangle_J,
\]
the Gram operator of $[\cdot,\cdot]$ with respect to
$\langle\cdot,\cdot\rangle_J$ is $J$ and $J$ is boundedly invertible. We
now prove part $(b)$ of the statement.
By Theorem \ref{8.8.7} there is only one Hilbert majorant, thus if there
are two positive inner products $\langle\cdot,\cdot\rangle_1$,
$\langle\cdot,\cdot\rangle_2$
whose associated norm define the Hilbert majorant, then the two norm must
be equivalent. Reasoning as in \cite{bognar}, the two Gram operators $G_j$
$j=1,2$ of $[\cdot,\cdot]$ with respect to $\langle\cdot,\cdot\rangle_j$,
$j=1,2$ are both boundedly invertible if and only if one of them is so.
Since we have previously shown that
$(b)$ holds for $G_1=J$ then $(b)$ holds for any other Gram operator.\\
Let us show the converse and assume that $(a)$ and $(b)$ hold. Then by
Theorem \ref{10.10.2}, $\mathcal K$ is decomposable and non-degenerate 
 thus, by Proposition \ref{4.4.1}, it admits a decomposition of the form
(\ref{krein}). By Proposition \ref{11.11.3}, $\mathcal K$ is a Krein space
if for every chosen decomposition the $J$-inner product makes $\mathcal K$
a Hilbert space or, equivalently, if $\tau_J$ coincides with $J$. First of
all we observe that since $G$ is boundedly invertible, by the closed
graph theorem we have that the
 Mackey topology coincides with $\tau$.   By Theorem \ref{9.9.7} we deduce
that $\tau$ is an admissible majorant and by Theorem \ref{9.9.5} $\tau$ is
also a minimal majorant and so $\tau\leq\tau_J$. However we know from
Proposition \ref{10.10.3} that $\tau_J\leq\tau$ and the conclusion
follows.
\end{proof}
\begin{Rk}
\label{10106rk}
{\rm Proposition \ref{10.10.6} says that in a Krein space all the
decomposition majorants are equivalent,
in other words,  all the $J$-norms are equivalent and will be called
natural norms
on $\mathcal K$. They define a Hilbert majorant called the strong topology
of $\mathcal K$.}
\end{Rk}
As an immediate consequence of the previous theorem we have:
\begin{Cy}
The strong topology of $\mathcal K$ equals the Mackey topology.
\end{Cy}
In the sequel we will always consider a Krein space $\mathcal K$ endowed
with the strong topology $\tau_M(\mathcal{K})$.
\\
\begin{Pn}
The strong topology $\tau_M(\mathcal{K})$ of the Krein space $\mathcal K$
is an admissible majorant.
\end{Pn}
\begin{proof}
By Proposition \ref{9.9.7} that the Mackey topology is
admissible and the fact that it is an admissible majorant is ensured by
(\ref{eqcs}). 
\end{proof}
\begin{Tm}\label{Tm:V3.4}
Let $\mathcal K$ be a quaternionic Krein space.  A subspace  $\mathcal{L}$
of  $\mathcal{K}$
is ortho-complemented if and only if it is closed and it is a Krein space
itself.
\end{Tm}
\begin{proof}
We assume that $\mathcal L$ is ortho-complemented. Then Corollary
\ref{cyIII2.5} shows that $\mathcal L$ is closed.
By Theorem \ref{Tm:V1.3}, $\mathcal K$ has a Hilbert majorant and thus we
can use the condition given in Proposition \ref{9.9.8} (b) to say whether
$\mathcal L$ is ortho-complemented. To this end, let us denote by
$G_{\mathcal{L}}$ the Gram operator defined by $[v,w]=\langle v,
G_{\mathcal{L}} w\rangle_J$, for $v,w\in\mathcal{L}$, where $J$ denotes
the fundamental symmetry of $\mathcal K$ associated with the chosen
decomposition. By Theorem \ref{Tm:V1.3}, the Gram operator $G$ is
boundedly invertible and thus, by Proposition \ref{9.9.8} (b) $\mathcal
L$ is ortho-complemented if and only if
${\rm Ran}(G_{\mathcal{L}}) =\mathcal{L}$ but, since $G_{\mathcal{L}}$ is
$J$-symmetric, this is equivalent to $G_{\mathcal{L}}$ boundedly
invertible and so, again by Theorem \ref{Tm:V1.3} to the fact that
$\mathcal L$ is a Krein space.
\end{proof}
Given a definite subspace $\mathcal{L}$ of a Krein space $\mathcal{K}$, it
is clear that the intrinsic topology $\tau_{int}(\mathcal{L})$ is weaker
than the topology induced by the strong topology $\tau_M(\mathcal{K})$
induces on $\mathcal{L}$. Thus we give the following definition:
\begin{Dn}
A subspace $\mathcal{L}$ of a Krein space $\mathcal{K}$ is said to be
uniformly positive (resp. negative) if $\mathcal L$ is positive definite
(resp. negative  definite) and
$\tau_{int}(\mathcal{L})=\tau_M(\mathcal{K})|\mathcal{L}$.
\end{Dn}
Note that the second condition amounts to require that $\mathcal{L}$ is
uniformly positive if  $[v,v]\geq c \|v \|_J^2$ for $v\in\mathcal{L}$
(resp.
$\mathcal{L}$ is uniformly negative if  $[v,v]\leq -c \|v \|_J^2$ for
$v\in\mathcal{L}$) where $c$ is a positive constant.
\begin{Tm}\label{T5.2}
Let $\mathcal K$ be a Krein space.
\begin{enumerate}
\item[(a)] A closed  definite subspace $\mathcal{L}$ of $\mathcal{K}$ is
intrinsically complete if and only if it is uniformly definite.
\item[(b)]  A semi-definite subspace $\mathcal{L}$ of $\mathcal{K}$ is
ortho-complemented if and only if it is closed and uniformly definite
(either positive or negative).
\end{enumerate}
\end{Tm}
\begin{proof}
The first statement follows from the fact that Proposition \ref{11.11.3} and the
closed graph theorem imply that a closed and definite subspace $\mathcal
L$ is intrinsically complete if and only if
$\tau_{int}(\mathcal{L})=\tau_M(\mathcal{K})|\mathcal{L}$ i.e. if and only
if $\mathcal L$ is uniformly definite.\\
By Proposition \ref{11.11.2} and Theorem \ref{Tm:V3.4}, a subspace
$\mathcal L$ is ortho-complemented if and only if  it is closed, definite
and intrinsically complete, i.e. if and only if $\mathcal L$ is uniformly
definite (either positive or negative). This completes the proof.
\end{proof}
\begin{Rk}
{\rm From the definition of uniformly positive (resp. negative) subspace, it follows that a subspace of $\mathcal K$ is uniformly positive (resp. negative) if so is its closure.}
\end{Rk}
Theorem \ref{T5.2} and the previous remark immediately give the following:
\begin{Cy}
A semi-definite subspace of $\mathcal K$ is uniformly definite if and only is its closure is ortho-complemented.
\end{Cy}
As a consequence of  Theorems \ref{Tm:V3.4} and \ref{T5.2} we also have the following
result, which was the main motivation for the present paper:
\begin{Tm}
\label{tm:sofsof}
Let $\mathcal K$ denote a quaternionic Krein space, and let
${\mathcal M}$ be a closed uniformly positive subspace of
$\mathcal K$. Then, ${\mathcal M}$ is a Hilbert space
and is ortho-complemented in $\mathcal K$: One can write
\[
\mathcal K= {\mathcal M}{[\oplus]} {\mathcal
M}^{[\perp]},
\]
and ${\mathcal M}^{[\perp]}$ is a Krein subspace of
$\mathcal K$.
\end{Tm}

\begin{proof}
The space is a Hilbert space by $(a)$ of Theorem \ref{T5.2}. That it is
orthocomplemented follows then from
Theorem \ref{Tm:V3.4}.
\end{proof}

\bibliographystyle{plain}


\begin{thebibliography}{10}

\bibitem{MR1333599}
S.L. Adler.
\newblock {\em Quaternionic quantum mechanics and quantum fields},
volume~88 of
  {\em International Series of Monographs on Physics}.
\newblock The Clarendon Press Oxford University Press, New York, 1995.

\bibitem{abcs2-milan}
D.~Alpay, V.~Bolotnikov, F.~Colombo, and I.~Sabadini.
\newblock Interpolation of {S}chur multipliers for slice hyperholomorphic
  functions.
\newblock In preparation.

\bibitem{MR2242728}
D.~Alpay, A.~C.~M. Ran, and L.~Rodman.
\newblock Basic classes of matrices with respect to quaternionic
indefinite
  inner product spaces.
\newblock {\em Linear Algebra Appl.}, 416(2-3):242--269, 2006.

\bibitem{as3}
D.~Alpay and M.~Shapiro.
\newblock Reproducing kernel quaternionic {P}ontryagin spaces.
\newblock {\em {Integral Equations and Operator Theory}}, 50:431--476,
2004.

\bibitem{azih}
T.~Ya. Azizov and I.S. Iohvidov.
\newblock {\em Foundations of the theory of linear operators in spaces
with
  indefinite metric}.
\newblock Nauka, Moscow, 1986.
\newblock {\rm (Russian). English translation: {\it Linear operators in
spaces
  with an indefinite metric}. John Wiley, New York, 1989}.

\bibitem{MR2394102}
J.A.~Ball, V.~Bolotnikov, and Q.~Fang.
\newblock Schur-class multipliers on the {A}rveson space: de
  {B}ranges-{R}ovnyak reproducing kernel spaces and commutative
  transfer-function realizations.
\newblock {\em J. Math. Anal. Appl.}, 341(1):519--539, 2008.

\bibitem{ball-fang}
J.A.~Ball and Q.~Fang.
\newblock {Nevanlinna-Pick interpolation via graph spaces and Krein-space
  geometry: a survey}.
\newblock {Preprint 2012. To appear in OTAA}.

\bibitem{MR756761}
J.A.~Ball and J.W.~Helton.
\newblock Beurling-{L}ax representations using classical {L}ie groups with
many
  applications. {II}. {${\rm GL}(n,\,{\bf C})$} and {W}iener-{H}opf
  factorization.
\newblock {\em Integral Equations Operator Theory}, 7(3):291--309, 1984.

\bibitem{bh-1986}
J.A.~Ball and J.W.~Helton.
\newblock Interpolation theorems of {P}ick--{N}evanlinna and {L}oewner
types
  for meromorphic matrix functions: parametrisation of the set of all
  solutions.
\newblock {\em {Integral Equations and Operator Theory}}, 9:155--203,
1986.

\bibitem{bognar}
J.~Bogn{\'a}r.
\newblock {\em Indefinite inner product spaces}.
\newblock Ergebnisse der Mathematik und ihrer Grenzgebiete, Band 78.
  Springer--Verlag, Berlin, 1974.

\bibitem{MR43:2}
N.~Bourbaki.
\newblock {\em \'{E}l\'ements de math\'ematique. {A}lg\`ebre. {C}hapitres
1 \`a
  3}.
\newblock Hermann, Paris, 1970.

\bibitem{MR83k:46003}
N.~Bourbaki.
\newblock {\em Espaces vectoriels topologiques. {C}hapitres 1 \`a 5}.
\newblock Masson, Paris, new edition, 1981.
\newblock \'El\'ements de math\'ematique. [Elements of mathematics].

\bibitem{bds}
F.~Brackx, R.~Delanghe, and F.~Sommen.
\newblock {\em Clifford analysis}, volume~76.
\newblock Pitman research notes, 1982.

\bibitem{MR2752913}
F. Colombo, I. Sabadini, and D.C. Struppa.
\newblock {\em Noncommutative functional calculus. Theory and
applications of slice hyperholomorphic functions}, volume 289 of
{\em Progress
  in Mathematics}.
\newblock Birkh\"auser/Springer Basel AG, Basel, 2011.


\bibitem{MR1364446}
M.A. Dritschel and J. Rovnyak.
\newblock Operators on indefinite inner product spaces.
\newblock In P.~Lancaster, editor, {\em Lectures on operator theory and
its
  applications ({W}aterloo, {ON}, 1994)}, volume~3 of {\em Fields Inst.
  Monogr.}, pages 141--232. Amer. Math. Soc., Providence, RI, 1996.

\bibitem{DS1}
N.~Dunford and J.~Schwartz.
\newblock {\em Linear operators}, volume~1.
\newblock Interscience, 1957.

\bibitem{MR0137500}
D. Finkelstein, J.M. Jauch, S. Schiminovich, and D. Speiser.
\newblock Foundations of quaternion quantum mechanics.
\newblock {\em J. Mathematical Phys.}, 3:207--220, 1962.

\bibitem{MR768240}
L.~P. Horwitz and L.~C. Biedenharn.
\newblock Quaternion quantum mechanics: second quantization and gauge
fields.
\newblock {\em Ann. Physics}, 157(2):432--488, 1984.

\bibitem{ikl}
I.S. Iohvidov, M.G. Kre\u{\i}n, and H.~Langer.
\newblock {\em Introduction to the spectral theory of operators in spaces
with
  an indefinite metric}.
\newblock Akademie--Verlag, Berlin, 1982.

\bibitem{MR2001c:15021}
M.~Karow.
\newblock Self-adjoint operators and pairs of {H}ermitian forms over the
  quaternions.
\newblock {\em Linear Algebra Appl.}, 299(1-3):101--117, 1999.

\bibitem{MR2433159}
L.~ Rodman.
\newblock Pairs of {H}ermitian and skew-{H}ermitian quaternionic matrices:
  canonical forms and their applications.
\newblock {\em Linear Algebra Appl.}, 429(5-6):981--1019, 2008.

\bibitem{MR2988215}
L.~ Rodman.
\newblock Stability of {I}nvariant {S}ubspaces of {Q}uaternion {M}atrices.
\newblock {\em Complex Anal. Oper. Theory}, 6(5):1069--1119, 2012.

\bibitem{Rud91}
W.~Rudin.
\newblock {\em Functional analysis}.
\newblock International Series in Pure and Applied Mathematics.
McGraw-Hill
  international editions, 1991.

\bibitem{schwartz}
L.~Schwartz.
\newblock Sous espaces hilbertiens d'espaces vectoriels topologiques et
noyaux
  associ\'{e}s (noyaux reproduisants).
\newblock {\em J. Analyse Math.}, 13:115--256, 1964.

\bibitem{MR2372587}
V.~Sergeichuk.
\newblock Canonical matrices of isometric operators on indefinite inner
product
  spaces.
\newblock {\em Linear Algebra Appl.}, 428(1):154--192, 2008.

\bibitem{viswanath_thesis}
K.~Viswanath.
\newblock {\em {Contributions to linear quaternionic analysis}}.
\newblock PhD thesis, {Indian Statistical Institute, Calcutta, India},
1968.

\bibitem{MR44:2067}
K.~Viswanath.
\newblock Normal operators on quaternionic {H}ilbert spaces.
\newblock {\em Trans. Amer. Math. Soc.}, 162:337--350, 1971.

\bibitem{MR97h:15020}
F.~Zhang.
\newblock Quaternions and matrices of quaternions.
\newblock {\em Linear Algebra Appl.}, 251:21--57, 1997.

\end{thebibliography}
\end{document}